\documentclass[11pt,a4paper]{article}
\usepackage[utf8]{inputenc}
\usepackage[T1]{fontenc}
\usepackage[english]{babel}
\usepackage{amsfonts,amssymb,epsfig}
\usepackage{graphicx,graphics,psfrag}
\usepackage{amsmath,amstext,amscd}
\usepackage{yhmath}

\usepackage{amsthm, cases}
\usepackage[percent]{overpic}
\usepackage[titletoc, title]{appendix}
\usepackage{enumitem}
\usepackage{empheq}
\usepackage{caption}
\usepackage{subcaption}

\usepackage{nicefrac}

\usepackage[margin=2cm]{geometry}

\usepackage{authblk}

\usepackage{multicol} 
\usepackage{setspace}

\usepackage[linesnumbered,ruled]{algorithm2e}

\usepackage{colortbl}

\usepackage{esint}
\usepackage{mathrsfs}
 
\usepackage{varwidth}
\usepackage{bbm} 
\usepackage{cite}

\usepackage{array}

\usepackage{dsfont}
\usepackage{ushort}
\usepackage{accents}
\usepackage{booktabs}

\usepackage{mathtools}
\mathtoolsset{showonlyrefs}

\numberwithin{equation}{section}

\usepackage{pgfplots}
\pgfplotsset{compat = newest}

\usepackage{color}

\usepackage[colorlinks=true,linkcolor=blue,citecolor=blue,urlcolor=blue,breaklinks]{hyperref}
\usepackage{url}
\newcommand{\lp}{\left(}
\newcommand{\rp}{\right)}
\newcommand{\eps}{\varepsilon}

\def\paragraph#1{{\bf #1\ }}

\numberwithin{equation}{section}

\newtheorem{remark}{Remark}[section]

\newcommand{\be}{\begin{equation}}
\newcommand{\ee}{\end{equation}}

\def\p{\partial}

\def\R{\mathbb{R}}

\def\1{\mathds{1}}

\def\d{\,\mathrm{d}}
\def\dx{\,\mathrm{d}x}
\def\dy{\,\mathrm{d}y}
\def\p{\,\partial}

\newcommand{\Id}{\mbox{Id}}
\newcommand{\cF}{\mathcal{F}}
\newcommand{\hcF}{\hat{\mathcal{F}}}
\newcommand{\cE}{\mathcal{E}}

\title{A Spectral-Based Method for Network-Formation PDEs}
\author[1]{Pedro Aceves-S\'{a}nchez}
\author[2]{Pierre Degond}
\author[3]{Sara Merino-Aceituno}
\author[4]{Claudia Wytrzens}

\affil[1]{Department of Mathematics, University of Arizona, USA\\
pedroas@arizona.edu}
\affil[2]{Institut de Math\'ematiques de Toulouse, CNRS \& Universit\'e
de Toulouse, Toulouse, France, pierre.degond@gmail.com}
\affil[3,4]{Faculty of Mathematics, University of Vienna, Oskar-Morgenstern-Platz 1, 1090 Vienna\\
sara.merino@univie.ac.at, claudia.wytrzens@univie.ac.at}

\begin{document}

\maketitle

\begin{abstract}
We propose and study a simple and scalable Fourier-based spectral method for a continuum model of network formation under periodic boundary conditions. The model provides the evolution of the pressure $p$ and the conductivity $m$ over time. The evolution of $p$ is given by an anisotropic Poisson equation, while the equation for $m$ contains three terms corresponding to a diffusion and an activation term of the network -- that depends on the gradient of the pressure -- as well as a relaxation term that acts as a decaying term. This system arises as a formal $L^2$-gradient flow of a non-convex energy functional. Our algorithm combines two ingredients: (i) a splitting method for the equation for $m$, where the activation and relaxation parts are solved analytically, and the diffusion part is solved via Fast Fourier Transform (FFT), and (ii) an FFT combined with the Conjugate Gradient (CG) method applied to the equation for the pressure. This makes the scheme easy to implement compared to implicit schemes and naturally extensible to three dimensions on uniform periodic grids. To showcase the method, we recover the previously documented influence of the activation strength $c$, the diffusion coefficient $D$, and the metabolic exponent $\gamma$ on the morphology of emergent networks, and report grid convergence results.\\

    
\end{abstract}

\tableofcontents

\section{Introduction}

\subsection{PDE System for Network Formation}

Biological transportation networks are characterised by their ability to efficiently balance cost, efficiency and resilience \cite{Corson2010intro, Barthelemy2011spatialnetworks_intro}. Such networks are fundamental processes and occur in various contexts, such as leaf venation in plants \cite{Corson2010intro, Malinowski2013leaf_intro}, blood circulation in mammals \cite{Sedmera2011blood_intro}, and the transport of electric charge in neural networks \cite{Michel1995neural_intro,Eichmann2005neural_intro}. Therefore, understanding these networks and the underlying mechanisms of their evolution is of great interest. Consequently, the aim of comprehending the emergence, formation, and adaptation of such networks has inspired many mathematical approaches to model these phenomena, ranging from discrete graph-based descriptions \cite{Hu2013discrete, Barabasi1999discrete_intro, Barthelemy2011spatialnetworks_intro} to continuum models based on partial differential equations (PDEs) \cite{Hu2019tensor}.

In this paper, we focus on a PDE model describing network formation in porous media first introduced in \cite{Hu2013_lecture_notes}. Particularly, we will propose a numerical method that is simpler compared to those proposed previously in the literature and that has the potential to be extended to three dimensions. The model describes the change of pressure $p=p(t,x,y) \in \R$ and conductivity $m=m(t,x,y)\in \R^2$ over the space variables $x,y \in \R$ and time $t \in [0,\infty)$. We have
\begin{subequations} \label{eq:networkSystem}
\be
-\nabla \cdot \left[  \left( (m\otimes m) + r \Id \right) \nabla p\right] =S, \label{eq:p}
\ee
\begin{equation}
\p_t m - D^2 \Delta m -c^2 (m \cdot \nabla p) \nabla p  + |m|^{2 (\gamma -1)}m =0, 
\label{eq:m}
\end{equation}
\end{subequations}
where the tensor product $\otimes$ is defined as $m\otimes m := \left( m_i m_j\right)_{i,j\in \{1,2\}}$,   $r(x,y)\geq r_0 >0$ is the isotropic background permeability of the medium,  $\gamma \in \R$ is the metabolic exponent,  $D>0$ is the diffusive constant and $c>0$ is the activation parameter. Furthermore, $\Id$ denotes the $(2\times2)$- identity matrix and $S$ represents the sources and sinks in the medium. In the sequel, we assume periodic boundary conditions for $p$ and $m$ and time-independent $S=S(x,y)$, with total integral $\int S(x,y) \d x \d y=0$ (where the integral is in the whole domain) since it is a solvability condition for \eqref{eq:p}. Specifically, for the numerical experiments (Section \ref{sec:results}), we will consider a rectangular domain $\Omega$ with periodic boundary conditions.

\bigskip

\paragraph{Meaning of the equations.}
In the elliptic equation for the pressure $p$ \eqref{eq:p}, the term $(m~\otimes~m~+~ r~\Id)$ is a permeability tensor depending on space $(x,y)$ and time $t$. This elliptic equation actually arises by combining the quasi-incompressibility of the fluid with Darcy's law. Darcy's law links the velocity $u$ and the pressure $p$ of a fluid in a porous medium in the form
\be
u(t,x,y)=-K(t,x,y)\nabla p(t,x,y), \label{eq:darcy}
\ee
where $K(t,x,y)$ is the conductivity of the medium, in our case $K(t,x,y)=(m\otimes m + r\Id)$. Now, assuming the quasi-incompressibility of the fluid and the conservation of mass, we have
\be 
\nabla \cdot u(t,x,y)=S(x,y).
\label{eq:masscons_incomp}
\ee
Substituting Darcy's law \eqref{eq:darcy} into relation \eqref{eq:masscons_incomp} leads to the equation for the pressure \eqref{eq:p}.

Additionally, note that the principal directions of the conductivity tensor $(m\otimes m + r\Id)$ are given by $m/|m|$ and $m^\perp$ with corresponding eigenvalues $r(x)+|m|^2$ and $r(x)$, respectively. This can also be interpreted as an increase of the principal permeability by $|m|^2$ in the direction of the conductivity $m$ while the flow only `feels' the background permeability $r(x)$ in the direction of $m^\perp$. 

Furthermore, equation \eqref{eq:m} is a parabolic reaction-diffusion equation which governs the evolution of the network conductivity $m$. It consists of the diffusive term $D^2\Delta m$ that accounts for random effects in the network; the activation term $c^2(m \cdot \nabla p) \nabla p$ that acts as a pumping term and thus the driving force to extend the network in the direction of the pressure gradient $\nabla p$; lastly, the algebraic relaxation term $|m|^{2(\gamma-1)}m$ acts as a decaying term and represents the functional derivative of the metabolic cost required to sustain the network.

The existence of weak solutions for \eqref{eq:networkSystem} has been established 
under various assumptions in a series of works 
\cite{Haskovec2016,Albi2017,Haskovec2014,Albi2016}, all within the 
regime $\gamma \geq \nicefrac{1}{2}$ and $r=r(x)>r_0>0$. The long-term behaviour of these 
solutions and their stability with respect to the diffusivity constant $D$ have 
been further investigated in \cite{Haskovec2014,Haskovec2016,Albi2017}; in 
particular, for small $D$ the zero steady-state is Turing unstable, giving rise 
to network patterns.

\medskip
The system for the pressure $p$ and the conductivity $m$ \eqref{eq:p}--\eqref{eq:m} is coupled, non-local, and it can be interpreted as the formal $L^2$-gradient flow induced by the energy functional $\cE[m]$, given by 
\begin{equation}
    \cE[m] := \frac{1}{2} \int_\Omega \left( D^2 |\nabla m|^2 + \frac{|m|^{2\gamma}}{\gamma} +c^2 (m \cdot \nabla p[m])^2 +c^2 r(x)|\nabla p[m]|^2 \right) \d x,   
    \label{eq:energy}
\end{equation}
which is highly non-convex \cite{Haskovec2014, Haskovec2016}. Here, $p[m] \in H^1(\Omega)$ denotes the unique solution of the Poisson equation \eqref{eq:p} corresponding to a given conductivity $m$, under periodic boundary conditions. The first and second terms in \eqref{eq:energy} represent the diffusion and relaxation energies, respectively, while the last two terms describe the energy associated with the interaction between the network and the fluid flow. In \cite{Haskovec2014}, the authors proved that the energy functional is non-increasing along smooth solutions and that for weak solutions, the following energy dissipation inequality holds:
$$\cE[m(t)]\leq \cE[m(0)] \qquad \text{for all } t \geq 0.$$

The system \eqref{eq:networkSystem} was derived from a discrete-graph model introduced by Hu and Cai \cite{Hu2013discrete}, which describes the evolution of conductivities on graphs subject to Kirchhoff's law and energy minimisation principles. Both approaches, the discrete (e.g., \cite{Hu2013discrete,Haskovec2019Auxin, Haskovec2019Murray}) and the continuous (e.g., \cite{Albi2016, Albi2017, Astuto2022comparison, Haskovec2014,Hu2019tensor}) one, as well as connections between these  (e.g., \cite{Burger2018,Haskovec2018ODE,Haskovec2019Rigorouslimit}), have been extensively studied resulting in a variety of analytical results on e.g., existence of solutions, stability and steady states, a rigorous continuum limit of the discrete model as well as numerical investigations illustrating arising network patterns.

\subsection{Contributions and Scope}

We propose a Fourier-based spectral algorithm for system \eqref{eq:networkSystem} on periodic domains that is:
(i) simple to implement (no assembly or inversion of large sparse matrices), (ii) high-order accurate in space, and (iii) naturally extensible to 3D uniform grids. The third point is particularly important when it comes to applications, as many biological networks, such as vascularisation, occur in three dimensions.

 Across the literature, system \eqref{eq:networkSystem} has been discretised using mixed/standard finite elements and finite differences, paired with time integrators tailored to stiffness: mixed finite elements in space with an implicit–explicit first‑order Euler scheme and adaptive time stepping (and Dirichlet/Neumann boundary conditions) in \cite{Haskovec2016}; finite elements with Crank–Nicolson for diffusion, explicit treatment of the remaining terms, and very small time steps to ensure energy decay in \cite{Albi2016}; operator splitting for the conductivity equation combining explicit Euler (activation), implicit Euler with bisection (relaxation), and exact diffusion, plus a preconditioned GMRES solver for the pressure in \cite{Albi2017}; fully implicit second‑order and semi‑implicit first‑order energy‑decaying schemes with a Conjugate Gradient solver for the Poisson equation in \cite{Fang2019implicit}; ADI (Alternating Direction Implicit method) and symmetric ADI schemes on uniform Cartesian grids for the tensor model in \cite{Astuto2022comparison,Astuto2025Asymmetry}; and low‑order finite elements with semi‑implicit backward Euler time stepping and Richardson extrapolation with accuracy assessment in \cite{Astuto2024FiniteElement}.  

This extensive and varied literature is due to the challenges in simulating system \eqref{eq:networkSystem}:  the non-linear nature of the system; the degeneracy for $D=0$; the stiffness and ill-conditioning issues -- for larger values of the background conductivity $r$ the stiffness increases, but as $r\to 0$ the problem becomes ill-conditioned (particularly, for fully implicit 2D solvers). Moreover, it has been reported that the solutions are very sensitive to the parameters $\gamma,D,r$ and spatial mesh sizes. Our method is not unconditionally energy-stable and inherits a time-step restriction from the activation subdynamics. Nonetheless, in the periodic setting, it delivers accurate and scalable simulations using FFTs, providing a practical computational tool that can complement implicit and energy-decaying schemes.

\subsection{Structure of this Paper}

In Section \ref{sec:numerical_method}, details of the proposed numerical method are given. Section \ref{sec:results} presents numerical experiments and parameter studies, and Section \ref{sec:gridconv} reports grid convergence results. Section \ref{sec:conclusion_network_formation} provides conclusions and outlines directions for future work.

\section{Fourier--Based Spectral Numerical Method} \label{sec:numerical_method}

\subsection{Key Idea and Motivation}
Our algorithm, similar to the one presented in \cite{Albi2017}, uses a splitting method for the time discretisation of equation \eqref{eq:m}. Specifically, at each time step we compute the conductivity $m$ by updating separately the values for the activation, relaxation and the diffusion terms. The key of our method, however, is in the resolution of the equation for $p$ \eqref{eq:p}. The goal is to use for this a Fourier-based spectral method -- implemented via the Fast Fourier Transform (FFT). This is motivated by computational and theoretical considerations. Since we formulate the model with periodic boundary conditions, implementation via FFT is a rather natural choice as it provides an efficient way to compute the spatial derivatives for $p$ and $m$ and to solve the elliptic Poisson equation for the pressure $p$, which reduces the computational cost from $\mathcal{O} (N^2)$ to $\mathcal{O}(N \log N)$ for $N$ grid points. Hence, a huge advantage of applying FFT compared to Finite Element Methods (FEM) is that the differential operators can be encoded in Fourier space -- leading to multiplications in this space -- and thus avoids the assembly and inversion of large sparse matrices as in FEM. Moreover, spectral methods such as FFT achieve higher accuracy for smooth enough solutions on coarser meshes compared to FEM. 
Furthermore, in the setup of uniform meshes and periodic boundary conditions, an extension to higher dimensions will be less complex and thus easier via FFT than via FEM.

\subsection{Discretisations and Notations}

Let $(t^k)_{k=0}^n$ with $n\geq1$ and $0\leq t^0 < ... <t^k<t^{k+1} < ... <t^n \geq T$ be a discretisation of the time interval $[0,T]$, with $\Delta t^k = t^{k+1}-t^k$.

Given a bounded rectangular domain $\Omega = [0,L_x]\times [0,L_y]$, for $L_x,\ L_y>0$, we discretise the domain via a uniform rectangular mesh $\mathcal{T}_h$ with $h=(h_x,h_y) \in \R^2$ and $h_x, h_y>0$ being the mesh sizes in $x$ and $y$ direction. Thus, we define the associated set of vertices to the given mesh $\mathcal{T}_{h}$ as $\{(x_i,y_j): 0\leq i \leq N_x, 0 \leq j \leq N_y\}$, where $N_x+1$ and $N_y+1$ denotes the number of elements along the $x$-axis and $y$-axis respectively. If we consider given mesh sizes $h_x$ and $h_y$ of the partition $\mathcal{T}_{h}$, we have $h_x=\Delta x= \frac{L_x}{N_x}$ and $h_y=\Delta y= \frac{L_y}{N_y}$, we obtain $0=x_0 < x_1 < ... <x_i<...<x_{N_x}=L_x$ with $x_{i+1}=x_i+\Delta x$ as partition of $[0,L_x]$ and analogously $y_{j+1}=y_j+ \Delta y $ with $y_0=0$ and $y_{N_y}=L_y$ as partition of $[0,L_y]$.

We denote by $p_h^k$ the time and space discretised solution of \eqref{eq:p} and by $m_h^k$ the time and space discretised solution of \eqref{eq:m}. We denote by $p_{ij}^k$ the approximated discretised solution of the pressure $p$ at node $(x_i,y_j)$ at time $t^k$, i.e., $p_{ij}^k \approx p(t^k, x_i,y_j)$ and analogously for the conductivity $m$, i.e., $m_{ij}^k\approx m(t^k,x_i,y_j)$.

\subsection{Splitting for the Conductivity Equation}
We split \eqref{eq:m} into activation, relaxation, and diffusion substeps.
\subsubsection{Activation Term}
The activation part of the parabolic equation of $m$ \eqref{eq:m} is given by 
\be
\p_t m = c^2 \lp m \cdot \nabla p \rp \nabla p = c^2 \lp \nabla p \otimes \nabla p \rp m.
\ee
For $|\nabla p(t,x,y)|\neq 0$ the exact solution reads 
\be
m(t,x,y)= \frac{1}{|\nabla p|^2}\left[ \lp (\nabla p)^\perp \otimes (\nabla p)^\perp \rp + \lp \nabla p \otimes \nabla p \rp \exp(c^2|\nabla p|^2 t)\right] m_0(x,y),
\label{eq:sol_activation}
\ee
with $m_0(x,y)=m(0,x,y)$.
If $|\nabla p(t,x,y)|=0$, then $\partial_t m(t,x,y)=0$ and in this case $m(t,x,y)=m_0(x,y)$. \\

The time and space discretised version of \eqref{eq:sol_activation} reads
\be
m_{ij}^{k+1}=\begin{cases}
    m_{ij}^k \qquad \text{if} \quad
    |\nabla p_{ij}^k|^2=0, \\
\frac{m_{ij}^k}{|\nabla p_{ij}^k|^2} \left( \left( \nabla p_{ij}^{k}\right)^\perp \otimes \left(\nabla p_{ij}^{k}\right)^\perp +\left( \nabla p_{ij}^k
\otimes \nabla p_{ij}^k\right) \exp\left(c^2 |\nabla p_{ij}^k|^2 \Delta t^k \right) \right) \qquad \text{else},
\label{eq:m_activation_discrete}
\end{cases}
\ee
where $m_{ij}^k=m(t^k,x_i,y_j)$ and $p_{ij}^k=p(t^k,x_i,y_j)$ are the numerical approximations of $m$ and $p$ at time $t^k$ at node $(x_i,y_j)$ as described in the beginning of the section. 

\begin{remark} \label{rem:Delta_t}
    The term $\exp(c^2 |\nabla p|^2 \Delta t)$ in expression \eqref{eq:sol_activation} creates stiffness, which restricts the admissible time-step size $\Delta t$. In particular, $\Delta t$ must be chosen sufficiently small so that this term remains sufficiently bounded and can be computed reliably in numerical simulations.
\end{remark}

\subsubsection{Relaxation Term}
The metabolic cost term for the evolution of the conductivity $m$ is given by the equation
\be
\p_t m = -|m|^{2(\gamma-1)}m,
\ee
whose solution is
\be
m(t,x,y)=\begin{cases}
    m_0(x,y)\exp(-t) \qquad \text{ if} \quad \gamma=1, \\
    m_0(x,y) \qquad \qquad \qquad \text{else if} \quad|m_0|=0, \\
    \frac{1}{|m_0(x,y)|}\lp |m_0(x,y)|^{2(1-\gamma)} - 2 (1- \gamma) t \rp^{\nicefrac{1}{2(1-\gamma)}}m_0(x,y) \qquad \text{otherwise.}
\end{cases}
\label{eq:sol_relaxation}
\ee

We consider the following time and space discretisation for this solution: 
\be
m_{ij}^{k+1}=\left\{
\begin{array}{ll}
 m_{ij}^k \exp(-\Delta t^k) & \mbox{ if $\gamma=1$},\\
m_{ij}^{k} & \mbox{ else if $|m_{ij}^k|=0$},\\
a^k_{ij}\1_{b^k_{ij}\geq 0} & \mbox{ otherwise,}
\end{array}
\right.
\label{eq:m_relaxation_discrete}
\ee
where
\begin{eqnarray}
    a^k_{ij}&:=&\frac{m_{ij}^k}{|m_{ij}^k|}\left( |m_{ij}^k|^{2(1-\gamma)}-2(1-\gamma)\Delta t^k\right)^\frac{1}{2(1-\gamma)},\\
    b^k_{ij}&:=&|m_{ij}^k|^{2(1-\gamma)}-2(1-\gamma)\Delta t^k,
\end{eqnarray}
and $\1$ is the indicator function. 

The indicator function in the last case of \eqref{eq:m_relaxation_discrete} ensures that the solution stays real, as negative values for $b_{ij}^k$ might lead to complex solutions, depending on the value of $\gamma \in [0,1)$. This approximation is necessary since, if we wanted to stop the computation at the time when the solution hits zero, this would put the following constraint on the time step
$$\Delta t \leq \frac{|m_0(x)|^{2(1-\gamma)}}{2(1-\gamma)},$$
which, for $\gamma \in [0,1)$, leads to very restrictive time steps.

\subsubsection{Diffusion Term}
The diffusive term is given by
\be
\partial_t m = D^2\Delta m. \label{eq:diffusion}
\ee
Since we consider periodic boundary conditions, we can apply the Fourier transform $\mathcal{F}$. We denote the Fourier-space variables by $\xi_x \in \R$ and $\xi_y\in \R$. With these variables, the Laplacian $\Delta$ in Euclidean space transforms into $-[\xi_x^2+\xi_y^2]$ in Fourier space. Hence, we are able to transform equation \eqref{eq:diffusion} into 
\begin{equation}
    \frac{d}{dt}\cF(m)=D^2\left[-\xi_x^2-\xi_y^2\right] \cF(m).
\end{equation}
These equations can be solved explicitly, and thus we obtain as solution for $m$
\begin{equation}
     m(t)=\cF^{-1}\left(e^{D^2\left[-\xi_x^2-\xi_y^2\right]t}\frac{\cF(m)(t^0)}{e^{D^2\left[-\xi_x^2-\xi_y^2\right]t^0}} \right), \label{eq:m_diffusion_FT}
\end{equation}
where $\cF^{-1}$ denotes the inverse Fourier transform and $t^0$ is the initial time. 

Next, to approximate this solution for $m$ numerically, we will use the Fast Fourier Transform (FFT) as approximation of the Fourier transform. Thus, we will consider the space and time discretised version of equation \eqref{eq:m_diffusion_FT} in the following. Remember that we consider the space $[0,T]\times \Omega$ where $\Omega = [0,L_x]\times [0,L_y]$. We denote the discrete version of the Fourier transform as $\hcF$ and its inverse by $\hcF^{-1}$. Moreover, we define the discretised Fourier-space variables $\hat{\xi}_x$ and $\hat{\xi}_y$ such that $\hat{\xi}_{x}(\ell)$ for $0\leq \ell < N_x$ is given by 
 
\be
    \hat\xi_x(\ell):= \left\{ \begin{array}{ll}
    \frac{2\pi \ell}{L_x} & \mbox{ for $0\leq \ell \leq \frac{N_x}{2}$,} \\
    \frac{2\pi (\ell-N_x)}{L_x} & \mbox{ for $\frac{N_x}{2}<\ell<N_x$},
    \end{array}
    \right.
\ee
where $N_x$ is the number of nodes used for the FFT (i.e., in our case $N_x=\nicefrac{L_x}{\Delta x}$ and analogously for $\hat{\xi}_y(j)$ as $0\leq j < N_y$. Thus, the time and space discretised solution $m_{\ell j }^k$ to the diffusion equation \eqref{eq:diffusion} reads

\begin{equation}
m_{\ell j}^{k+1}=\hcF^{-1}\left(e^{D^2\left[-\xi_x(\ell)^2-\xi_y(j)^2\right]\Delta t^k}\cF(m_{\ell j}^k)\right). \label{eq:diffusion_space_discrete}
\end{equation}
For the implementation of the approximation of the diffusion term via FFT, we use the FFTW library developed by Frigo and Johnson \cite{FFTW} and apply it in two dimensions. For further details on the FFTW and its implementation, see \cite{FFTW}.

\medskip
\subsection{Equation for the Pressure \texorpdfstring{$p$}{p}: FFT+CG} \label{sec:Eq_pressure}

The method we use to compute the equation for the pressure $p$ \eqref{eq:p} combines the Fourier transformation and the Conjugate Gradient (CG) method. \\
As for the diffusive part of the equation for $m$, i.e., \eqref{eq:diffusion}, we  apply the FFT to transform equation \eqref{eq:p} into   
\be
- \mathcal{F}^{-1} \left[ i \left( \begin{array}{c}
     \xi_x \\
     \xi_y
\end{array} \right) \cdot \mathcal{F}\left( \left( m\otimes m +r \Id \right) \mathcal{F}^{-1}\left[ i \left( \begin{array}{c}
     \xi_x \\
     \xi_y
\end{array} \right) \mathcal{F}(p) \right]
\right)
\right] =  S, 
\ee
where $\mathcal{F}$ is the continuous Fourier transformation, $\mathcal{F}^{-1}$ the inverse Fourier transformation, $i$ is the imaginary unit, and $\xi_x$ and $\xi_y$ are the Fourier-space variables. \\

Setting $\hat{\mathcal{F}}$ and $\hat{\mathcal{F}}^{-1}$ as the discrete (inverse) Fourier transformation, the time and space discretised version of the above  equation reads
\be
- \hat{\mathcal{F}}^{-1} \left[ i \left( \begin{array}{c}
     \hat\xi_{x}(\ell) \\
     \hat\xi_y(j)
\end{array} \right) \cdot \hat{\mathcal{F}}\left( \left( m_{\ell j}^k\otimes m_{\ell j}^k +r \Id \right) \hat{\mathcal{F}}^{-1}\left[ i \left( \begin{array}{c}
     \hat\xi_x(\ell) \\
     \hat \xi_y(j)
\end{array} \right) \hat{\mathcal{F}}\left(p_{\ell j}^{k+1}\right) \right]
\right)
\right] =  S_{\ell j}. \label{eq:p_discrete}
\ee
Here, $S_{\ell j}$ is the numerical approximation of the source-and sink terms, i.e., $S_{\ell j} \approx S(x_\ell, x_j)$.\\

Now, in order to compute $p^{k+1}_{\ell j}$, we will apply the CG method to equation $\eqref{eq:p_discrete}$ with $p^k_{\ell j}$ as initial data for the iterative method. In general, the CG method solves symmetric, positive definite linear systems of the form $Ax=b$ with an iterative procedure, where $A \in \R^{n\times n}$ and $x,b \in \R^n$ \cite{BD2008}. Note that for the numerical computation of the CG method, it is not necessary that the matrix $A$ is explicitly available, since only the matrix-vector product will be used, and moreover, the CG algorithm terminates if the residual error is smaller than a given threshold \cite{BD2008}. We apply the CG method as described in Algorithm 10.1 in \cite{BD2008}. For further details on the CG method, we refer to Chapter 10.4 and, in particular, Algorithm 10.1 in \cite{BD2008}.

\subsection{Stopping Criteria and Energy Evaluation}

We impose two stopping criteria.
 One is based on the relative change of the conductivity $m$:
\be
\Delta m^k_h:=\frac{|| m^k_h-m^{k-1}_h||_2}{\Delta t^k ||m^{k-1}_h||_2 },
\ee
here $\|\cdot\|_2$ is the discrete $L^2-$norm.
We use this stopping criterion as a signal that the simulation is close to a stationary solution. 
The second stopping criterion terminates the algorithm if the discrete energy increases, that is, if $\cE^k_h\geq \cE^{k-1}_h$, with
\be
\cE_h^k(m^k_h) := \frac{D^2}{2}\|\nabla m_h^k\|_2^2 +\frac{r c^2}{2}\|\nabla p_h^k\|_2^2+ \frac{c^2}{2}\|m_h^k \cdot \nabla p_h^k \|^2_2+\frac{1}{2\gamma}\||m_h^k|^\gamma\|^2_2.
\ee

\subsection{Pseudocode}

Table \ref{alg:pseudocode} gives the pseudocode to simulate system \eqref{eq:networkSystem}.

\begin{algorithm}[ht!]
    \LinesNumbered
    \SetKwInput{Input}{Input}
    \SetKwInput{Output}{Output}
    \SetKwInput{Initialization}{Initalialization}
\SetKwBlock{while}{while}{}
    \Input{$\Omega=[0,L_x]\times[0,L_y]$, $r,c,D,\gamma$, $\Delta x,\Delta y$, $\Delta t$, initial $m_0$, source $S$ with $\int_\Omega S \dx \dy =0$, periodic BCs}
  
    \Output{$p^{k+1}_h,m^{k+1}_h$\\
    
    \text{ } \\
    }
    
     \while( $\Delta m_h^k:=\frac{||m_{h}^k-m_h^{k-1}||_{2}}{\Delta t^k ||m_h^{k-1}||_{2}}> \eps$ 
     and $\cE_h \lp m^{k+1}_h\rp$ $<$ $\cE_h\lp m^k_h \rp$ )
     {        \textbf{Step 1 (solve equation \eqref{eq:p}):} Compute $p_h^{k+1}$ with $S, m_h^{k}$ and $r$ given via solving \\
        $$ -\nabla \cdot [(rI+m_h^k\otimes m_h^k)\nabla p_h^{k+1}]=S$$ with FFT+CG as described in Section \ref{sec:Eq_pressure};\\
        \text{ }\newline
        \textbf{Step 2 (solve eq.\eqref{eq:m}
        via splitting algorithm):}\\

             \textit{\underline{Activation:}} (use explicit solution)\\ 
        compute gradient of $p$ via FFT \\ 
        compute solution of $m$ of activation term as stated in equation \eqref{eq:m_activation_discrete}

          \vspace{0.2cm}\textit{\underline{Relaxation:}} (use explicit solution) \\ 
          update $m_h^{k+1}$ according to equation \eqref{eq:m_relaxation_discrete}

          \vspace{0.2cm} \textit{\underline{Diffusion:}} (solved via FFT)\\
           compute $m_h^{k+1}$ via FFT and the inverse FFT as stated in equation \eqref{eq:diffusion_space_discrete}
         
		\text{ }\\
		
  }
	return solutions $p_h^{k+1},m_h^{k+1}$\\
{}
    \caption{FFT+CCG splitting algorithm for \eqref{eq:networkSystem}}
    \label{alg:pseudocode}
\end{algorithm}

\section{Numerical Experiments}\label{sec:results}
\subsection{Numerical Setup}

To showcase the numerical method presented in Section \ref{sec:numerical_method}, we will consider two types of experiments on $\Omega=[0,1]^2$ with periodic boundary conditions inspired by previously considered setups in the literature. Unless stated otherwise, we use $r=0.1$, $c=200$, a fixed time step $\Delta t=0.005$ conforming to Remark \ref{rem:Delta_t}, and a CG residual threshold $\|res^k\|\le 10^{-7}$. For the stopping criteria of the algorithm, we pick $\Delta m^k_h \leq 10^{-4}$. 
\medskip

We also consider three different initial conditions for the conductivity vector $m$:   
\begin{align}
m^A_{\mathrm{ini}}&=\big(\exp(-10^4(y-0.5)^2-50(x-0.5)^2),\,0\big),\label{eq:m_ini_normal}\\ 
m^B_{\mathrm{ini}}&=\big(\exp(-50(y-0.5)^2-50(x-0.5)^2),\,0\big),\label{eq:m_ini_symm}\\ %
m^{\mathrm{const}}_{\mathrm{ini}}&=(1,1).\label{eq:m_ini_const}
\end{align}

Note that all these initial conditions for the conductivity are periodic in the considered domain $\Omega$. Here $m_{ini}^{const}$ is isotropic, while $m_{ini}^A$ is anisotropic. The initial condition $m_{ini}^B$ is isotropic but not equal in its components. The maximum value of $m_{ini,1}^A$ and $m_{ini,1}^B$ in $\Omega$ is reached at the centre of the domain $\Omega$ at (0.5,0.5), and then it decays exponentially. However $m_{ini,1}^A$ decays faster along the $y$ axis. The initial data $m_{ini}^A$, $m_{ini}^B$ and $m_{ini}^{const}$ are inspired from the ones used in previous literature \cite{Albi2017, Astuto2022comparison}.

Inspired also from previous literature, we will consider two setups for the source $S$ in \eqref{eq:p}. In Section \ref{sec:1source4sinks}, we investigate experiments where we examine a squared domain with one source and four sinks. In Section \ref{sec:internalsink} we study the case where the whole domain acts as a sink.

\subsection{One Source and Four Sinks}\label{sec:1source4sinks}
In this setup, we consider one source and four sinks so that $S=S_4(x,y)$ reads 
\begin{equation}
   S_4(x,y)= S_{source}(x,y)-\sum_{i=1}^4 S_{sink_i},  \label{eq:S_4sinks}\\
\end{equation}
with
\begin{equation}
    S_{source}(x,y) = \exp(-1000((y-0.5)^2+(x-0.5)^2)), \label{eq:S_4sinks_source}\\
\end{equation}
and \\
 \begin{equation}
     S_{sink_i} =  \frac{1}{4} \exp(-1000((y-y_i)^2+(x-x_i)^2)) \qquad \text{for } i \leq 4. \label{eq:S_4sinks_sink}
 \end{equation}
Notice that indeed $\int_\Omega S(x,y)dxdy=0$ holds. Here, the source $S_{source}$ is centred at $(0.5, 0.5)$ and decays rapidly (exponentially) in an isotropic way. The sinks are located at different points, they also decay exponentially and are isotropic. Specifically, the sinks are located at $(x_1,y_1)=(0.2,0.2), (x_2,y_2)=(0.2,0.8), (x_3,y_3)=(0.8,0.2)$ and $(x_4,y_4)=(0.8,0.8)$, which results in a periodic setup, which is symmetric along the axes. For the simulations, we consider $c=200$ and $D=0.0025$, since these are values for which a clear network emerges, and an initial conductivity given by $m^A_{ini}$ in \eqref{eq:m_ini_normal}. This initial data will make the network want to develop towards the right and the left.

\begin{figure}
    \centering
    \begin{subfigure}[b]{0.48\textwidth}
         \centering
         \includegraphics[width=1.1\linewidth]{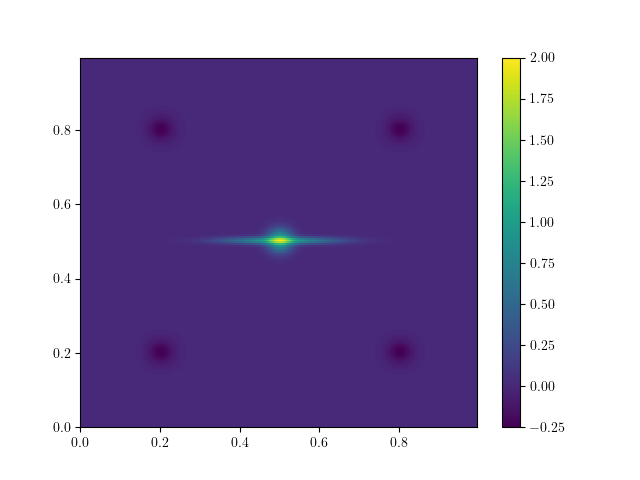}         
     \end{subfigure}
     \hfill
     \begin{subfigure}[b]{0.48\textwidth}
         \centering
         \includegraphics[width=1.1\linewidth]{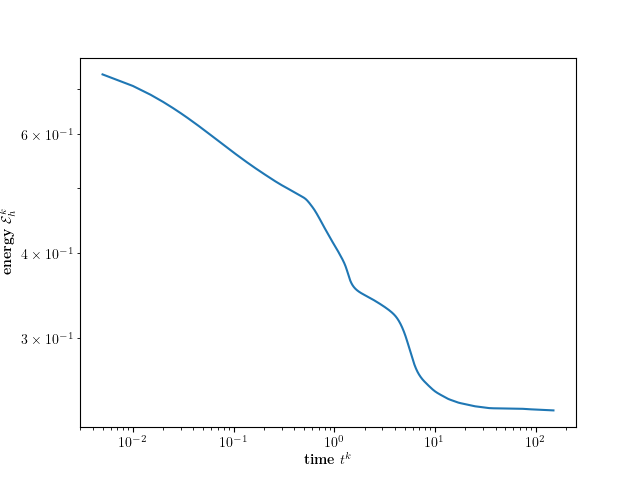}
    \end{subfigure}
       \caption{Left: Source and sink term $S_4(x,y)$ and initial conductivity $m_{ini}^A$. Right: Evolution of the discrete energy $\cE_h^k$ for $c=200, D=0.0025$ and $\gamma=1$.}
        \label{fig:4sinks_ini}
\end{figure}

Figure \ref{fig:4sinks_ini} depicts this initial setup for these simulations as well as the time evolution of the discretised energy $\cE_h^k$ for $\gamma=1$.

\begin{figure}[ht!]
    \hspace{-1.5cm}
     \centering
     \begin{minipage}{0.9\textwidth}
     \hfill
    \begin{subfigure}[b]{0.32\textwidth}
         \centering
         \includegraphics[width=0.9\linewidth,
                         trim=10mm 5mm 40mm 15mm, clip]         {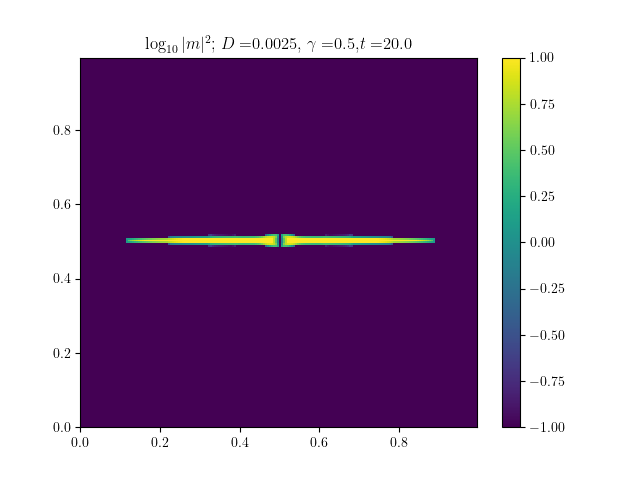}
         \caption{$\gamma=0.5$\\ \text{ }}
         \label{fig:4sinks_g05}
     \end{subfigure}
     \hfill
     \begin{subfigure}[b]{0.32\textwidth}
         \centering
         \includegraphics[width=0.9\linewidth,
                         trim=10mm 5mm 40mm 15mm, clip]         {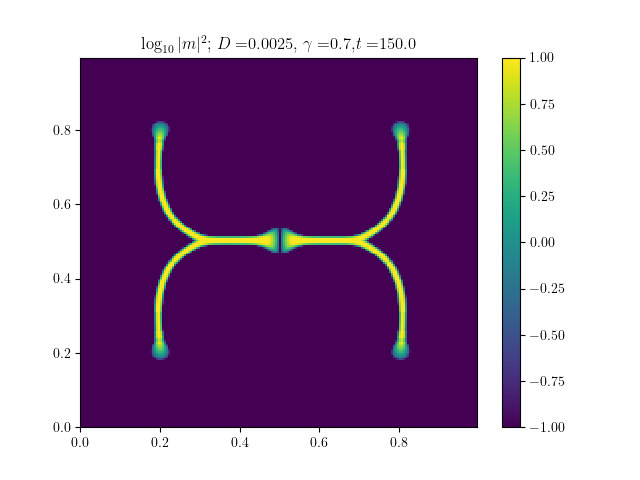}
         \caption{$\gamma=0.7$\\ \text{ } }
         \label{fig:4sinks_g07}
     \end{subfigure}
     \hfill
     \begin{subfigure}[b]{0.32\textwidth}
         \centering
         \includegraphics[width=0.9\linewidth,
                         trim=10mm 5mm 40mm 15mm, clip]
         {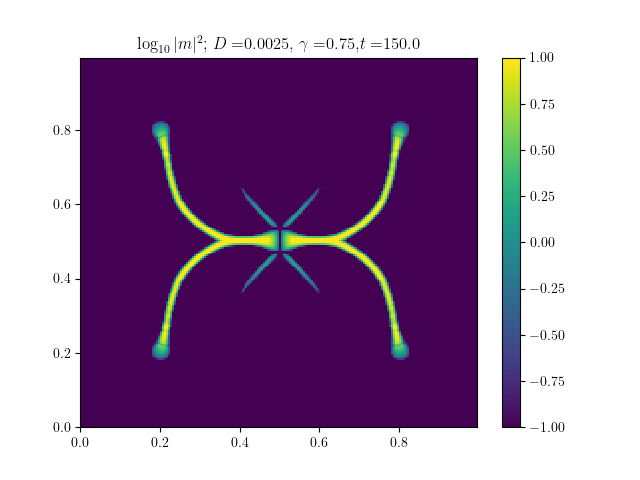}
         \caption{$\gamma=0.75$\\ \text{ } }
         \label{fig:4sinks_g075}
     \end{subfigure}
     \vspace{0.5cm}
     \hfill
     \begin{subfigure}[b]{0.32\textwidth}
         \centering
         \includegraphics[width=0.9\linewidth,
                         trim=10mm 5mm 40mm 15mm, clip]
         {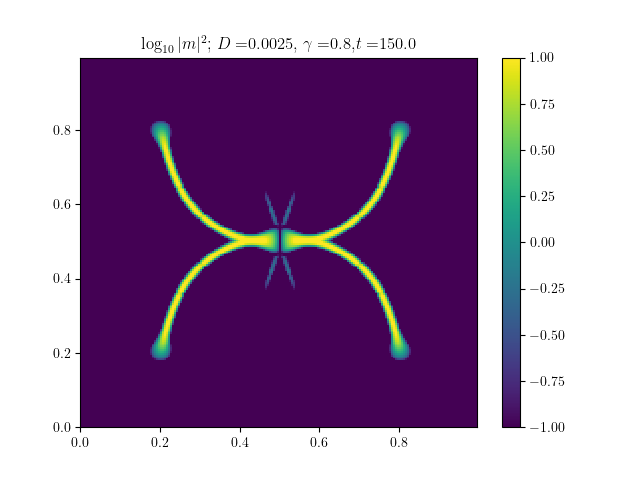}
         \caption{$\gamma=0.8$ }
         \label{fig:4sinks_g08}
     \end{subfigure}
\hfill
         \begin{subfigure}[b]{0.32\textwidth}
         \centering
         \includegraphics[width=0.9\linewidth,
                         trim=10mm 5mm 40mm 15mm, clip]
         {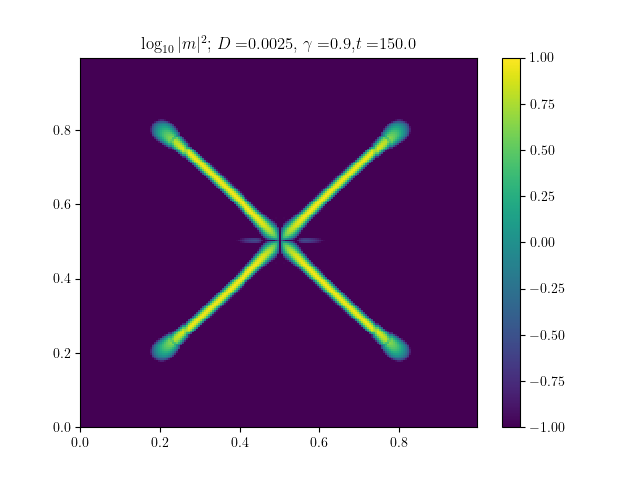}
         \caption{$\gamma=0.9$}
         \label{fig:4sinks_g09}
     \end{subfigure}
     \hfill 
    \begin{subfigure}[b]{0.32\textwidth}
         \centering
         \includegraphics[width=0.9\linewidth,
                         trim=10mm 5mm 40mm 15mm, clip]
         {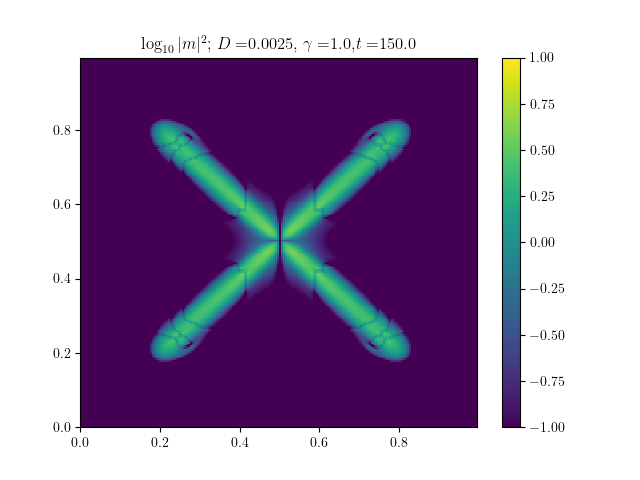}
         \caption{$\gamma=1$}
         \label{fig:4sinks_g1}
     \end{subfigure}
     \end{minipage}
       \begin{minipage}{0.05\textwidth}
         \begin{subfigure}[t]{1.5\textwidth}
         \centering
         \vspace{-1.4cm}
         \includegraphics[height=9cm]
         {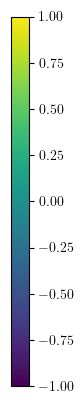}
     \end{subfigure}
     \end{minipage}

        \caption[Steady states of $\log_{10} |m|^2$ for $S_4$ and varying $\gamma$]{Steady states of $\log_{10} |m|^2$ for $c=200$, $D=0.0025$, $m^A_{ini}$, $S_4$ (one source, four sinks) and varying metabolic exponent $\gamma \in \{0.5,0.7,0.75,0.8, 0.9, 1\}$. We observe more symmetric networks and branches moving closer to the source as $\gamma$ increases.}
        \label{fig:4sinks_normal_m}
\end{figure}

\begin{figure}[ht!]
     \centering
    \begin{subfigure}[b]{0.32\textwidth}
         \centering
         \includegraphics[width=1.2\linewidth,
                         trim=0mm 0mm 0mm 0mm, clip]         {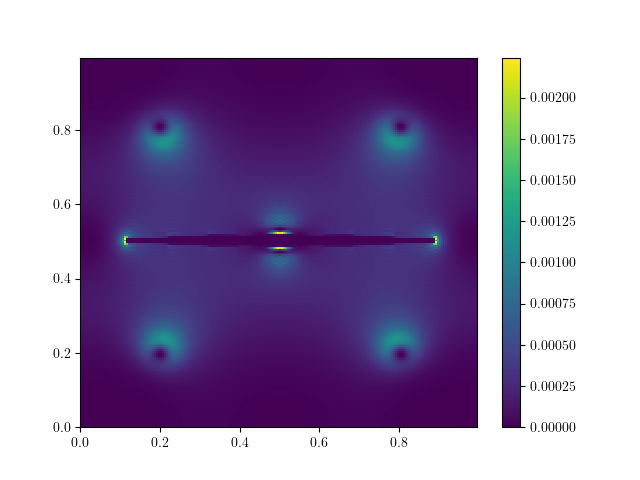}
         \caption{$\gamma=0.5$\\ \text{ }}
         \label{fig:4sinks_g05_gradp}
     \end{subfigure}
     \hfill
     \begin{subfigure}[b]{0.32\textwidth}
         \centering
         \includegraphics[width=1.2\linewidth,
                         trim=-3mm 0mm 0mm 0mm, clip]         {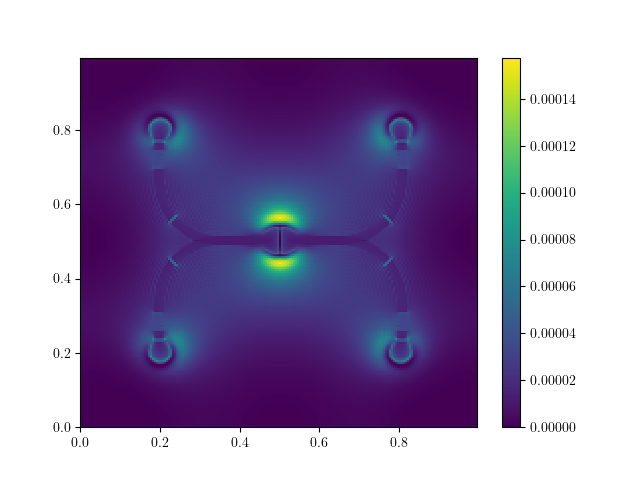}
         \caption{$\gamma=0.7$\\ \text{ }}
         \label{fig:4sinks_g07_gradp}
     \end{subfigure}
     \hfill
     \begin{subfigure}[b]{0.32\textwidth}
         \centering
         \includegraphics[width=1.2\linewidth,
                         trim=-3mm 0mm 0mm 0mm, clip]
         {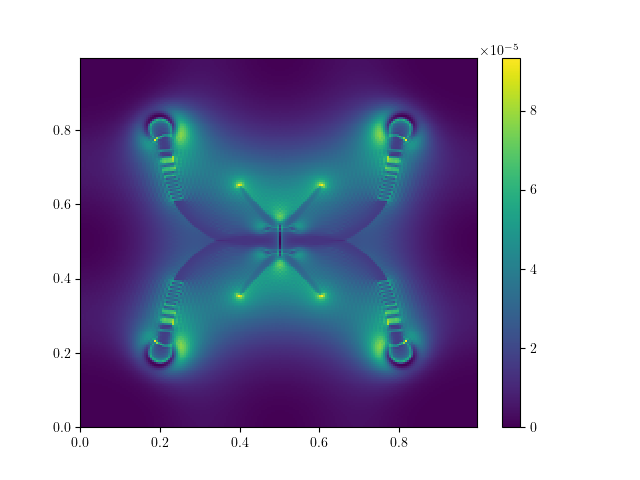}
         \caption{$\gamma=0.75$ \\ \text{ }}
         \label{fig:4sinks_g075_gradp}
     \end{subfigure}
     \hfill
     \begin{subfigure}[b]{0.32\textwidth}
         \centering
         \includegraphics[width=1.2\linewidth,
                         trim=0mm 0mm 0mm 0mm, clip]
         {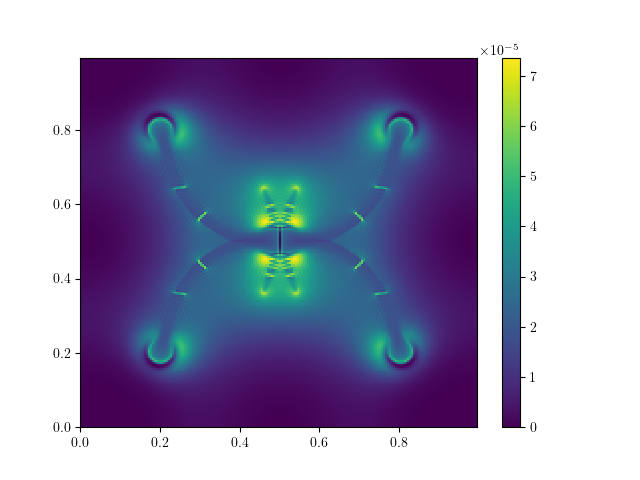}
         \caption{$\gamma=0.8$}
         \label{fig:4sinks_g08_gradp}
     \end{subfigure}
\hfill
         \begin{subfigure}[b]{0.32\textwidth}
         \centering
         \includegraphics[width=1.2\linewidth,
                         trim=-3mm 0mm 0mm 0mm, clip]
         {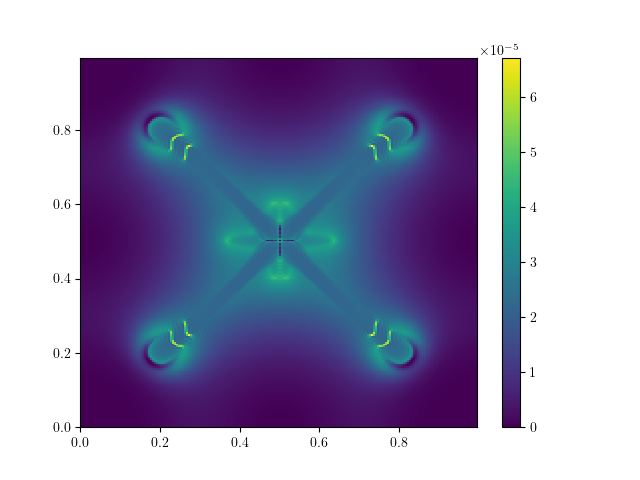}
         \caption{$\gamma=0.9$}
         \label{fig:4sinks_g09_gradp}
     \end{subfigure}
     \hfill 
    \begin{subfigure}[b]{0.32\textwidth}
         \centering
         \includegraphics[width=1.2\linewidth,
                         trim=-3mm 0mm 0mm 0mm, clip]
         {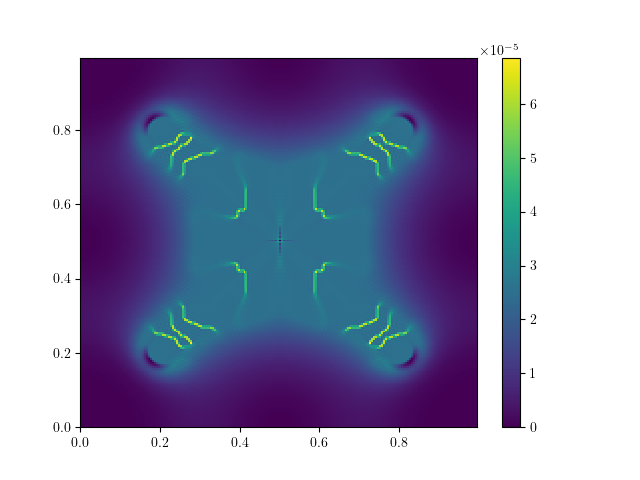}
         \caption{$\gamma=1$}
         \label{fig:4sinks_g1_gradp}
     \end{subfigure}
     
        \caption[Steady states of $|\nabla p|^2$ for $S_4$ and varying $\gamma$]{Steady states of $|\nabla p|^2$ for $c=200, D=0.0025$, $m_{ini}^A$, $S_4$ (one source, four sinks) and varying metabolic exponent $\gamma \in \{0.5,0.7,0.75,0.8, 0.9, 1\}$. }
        \label{fig:4sinks_normal_gradp}
\end{figure}

In Figure \ref{fig:4sinks_normal_m} and \ref{fig:4sinks_normal_gradp}, the corresponding simulation results for $\gamma \in \{0.5, 0.7, 0.75, 0.8, 0.9, 1\}$ for $\log_{10}|m|^2$ and for $|\nabla p|^2$ are presented. Note that we set $t_{end}=150$ and all simulations except for $\gamma=0.5$ stopped due to reaching $t_{end}$. For these simulations, the different energy evolutions over time look similar to the one depicted in Figure \ref{fig:4sinks_ini} for $\gamma=1$, and thus we assume that the results for $t_{end}=150$ are close to the steady states. For $\gamma=0.5$, the simulation broke at $t^k=20.545$ due to energy increase. This is the reason why we do not observe a network between the source and the sinks. We remind the reader that for general $r$, the existence of solutions has been proven only for $\gamma\geq\frac{1}{2}$. We observe that for $\gamma \in \{0.9, 1\}$ at the steady states, there are no branches but rather straight lines from the source to the sinks, and as $\gamma$ increases, the network is getting thicker. Regarding $\gamma \in \{0.7,0.75,0.8\} $, we notice that as $\gamma$ increases, the branching point of the network moves closer to the source in the centre. In general, we see that as $\gamma$ increases, the steady state becomes more diagonally symmetric, which means that the anisotropic effect of the initial conductivity $m^A_{ini}$ is lost.
The case $\gamma=1$ is a critical case. We see that in this case, the branches are significantly wider, indicating that diffusion becomes more dominant. 
In Figure $\ref{fig:4sinks_normal_gradp}$ we show the corresponding pressure fields. Notably, the maximum pressure decreases as $\gamma$ decreases. For larger $\gamma$ we observe localized segments within the branches where the pressure is higher.

\subsection{Internal Sink}\label{sec:internalsink}

For this section, we only consider an internal sink term, which means that the sink is constantly spread over the whole domain such that the total integral of the source and sink term is 0, i.e., $\int_{\Omega} S(x,y) \d x \d y = 0$. In the following the source and sink term $S=S_{int}(x,y)$ is given by 
\begin{align} \label{eq:SinternalSink}
S_{int}(x,y)=\exp(-1000((y-0.5)^2+(x-0.5)^2))-\frac{\pi}{1000}.
\end{align}

\begin{figure}[ht!]
\hspace{-1.5cm}
     \centering
     \begin{minipage}{0.9\textwidth}
     \hfill
    \begin{subfigure}[b]{0.32\textwidth}
         \centering
         \includegraphics[width=0.9\linewidth,
                         trim=10mm 5mm 40mm 15mm, clip]         {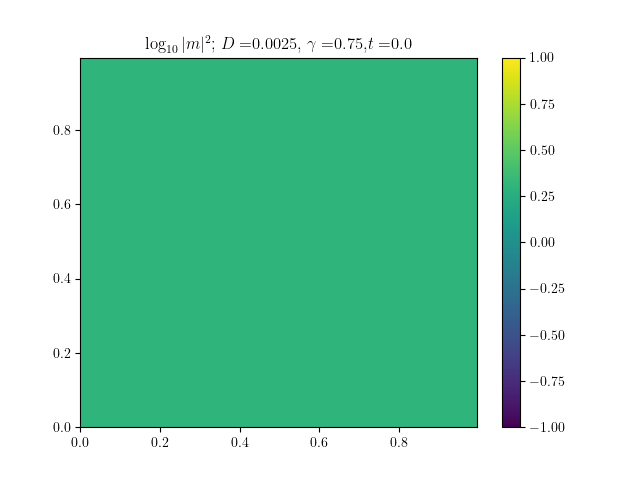}
         \caption{$t=0.0$ \\ \text{}}
         \label{fig:mini_t0}
     \end{subfigure}
     \hfill
     \begin{subfigure}[b]{0.32\textwidth}
         \centering
         \includegraphics[width=0.9\linewidth,
                         trim=10mm 5mm 40mm 15mm, clip]         {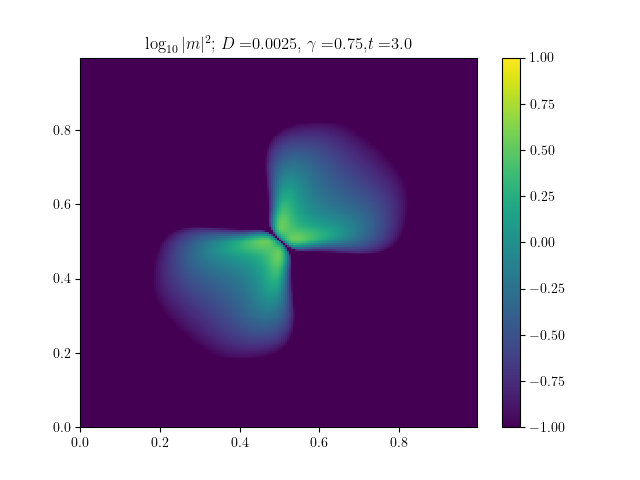}
         \caption{$t=3.0$ \\ \text{} }
         \label{fig:mini_t3}
     \end{subfigure}
     \hfill
     \begin{subfigure}[b]{0.32\textwidth}
         \centering
         \includegraphics[width=0.9\linewidth,
                         trim=10mm 5mm 40mm 15mm, clip]
         {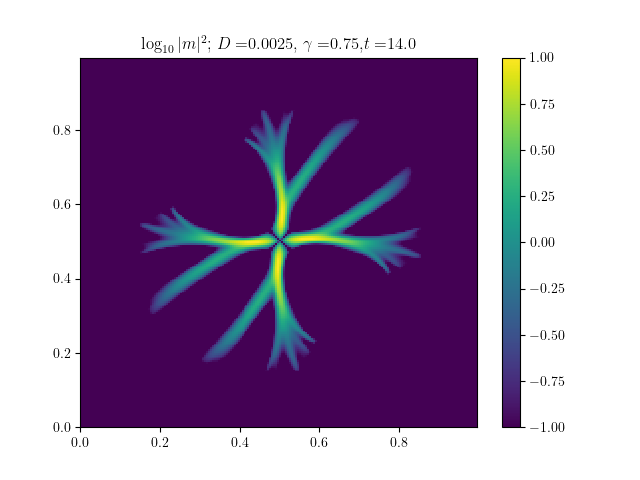}
         \caption{$t=14.0$ \\ \text{}}
         \label{fig:mini_t14}
     \end{subfigure}
     \hfill
     \begin{subfigure}[b]{0.32\textwidth}
         \centering
         \includegraphics[width=0.9\linewidth,
                         trim=10mm 5mm 40mm 15mm, clip]
         {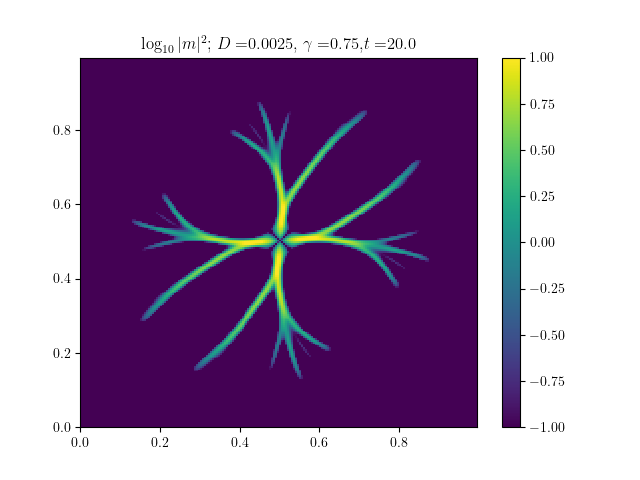}
         \caption{$t=20.0$}
         \label{fig:mini_t20}
     \end{subfigure}
\hfill
         \begin{subfigure}[b]{0.32\textwidth}
         \centering
         \includegraphics[width=0.9\linewidth,
                         trim=10mm 5mm 40mm 15mm, clip]
         {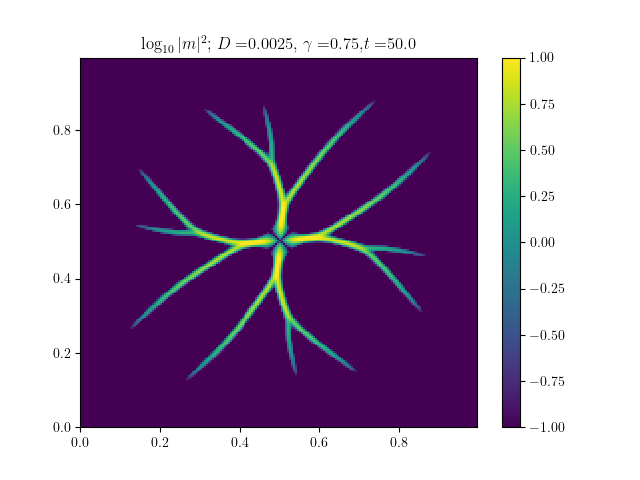}
         \caption{$t=50.0$}
         \label{fig:mini_t50}
     \end{subfigure}
     \hfill 
    \begin{subfigure}[b]{0.32\textwidth}
         \centering
         \includegraphics[width=0.9\linewidth,
                         trim=10mm 5mm 40mm 15mm, clip]
         {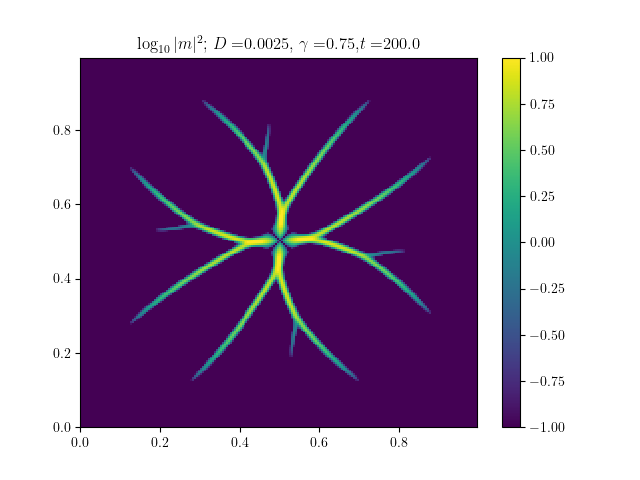}
         \caption{$t=200.0$}
         \label{fig:mini_t200}
     \end{subfigure}
     \end{minipage}
     \begin{minipage}{0.05\textwidth}
         \begin{subfigure}[t]{1.5\textwidth}
         \centering
         \vspace{-0.9cm}
         \includegraphics[height=9.1cm]
         {plots/internal_sink/constant_mini/c200_D00025_gamma075/skala.png}
     \end{subfigure}
     \end{minipage}
        \caption[Time evolution of network formation for constant initial conductivity $m$]{Time evolution of $\log_{10} |m|$ for $c=200, D=0.0025$ and $\gamma=0.75$, $m_{ini}^{const}$ and $S_{int}$ (internal sink). At times $t=3$ and $t=14$ a dense branching network is formed. But as time evolves, smaller branches vanish, and the structure simplifies such that only a few main branches remain at $t_{end}=200$.}
        \label{fig:mini_time_evolution}
\end{figure}

Figure \ref{fig:mini_time_evolution} shows the time evolution of the formation of the network for $c=200$,\linebreak $D=0.0025, \gamma=0.75$ starting from the constant initial condition $m_{ini}^{const}$ as stated in \eqref{eq:m_ini_const}. We observe that branches are evolving, but they vanish over time as the energy decreases, and only a few main branches remain at later times ($t=50$, $t=200$), indicating a transition towards a simpler, energetically favourable structure. Notice that, despite $m_{\mathrm{ini}}^{\mathrm{const}} = (1,1)$ being constant 
in space and $S_{\mathrm{int}}(x,y)$ being radially symmetric around $(0.5, 0.5)$, 
the resulting network does not exhibit reflection symmetry about either coordinate 
axis. Instead, symmetry is preserved only under reflection across the diagonals.

\begin{figure}[ht!]
     \centering
    \begin{subfigure}[b]{0.48\textwidth}
         \centering
         \includegraphics[width=1.1\linewidth]         {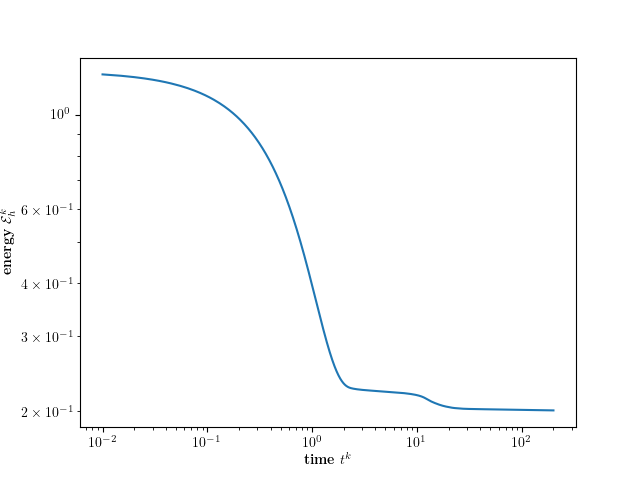}
         \caption{Energy evolution over time}
         \label{fig:mini_t0s}
     \end{subfigure}
     \hfill
     \begin{subfigure}[b]{0.48\textwidth}
         \centering
         \includegraphics[width=1.1\linewidth]         {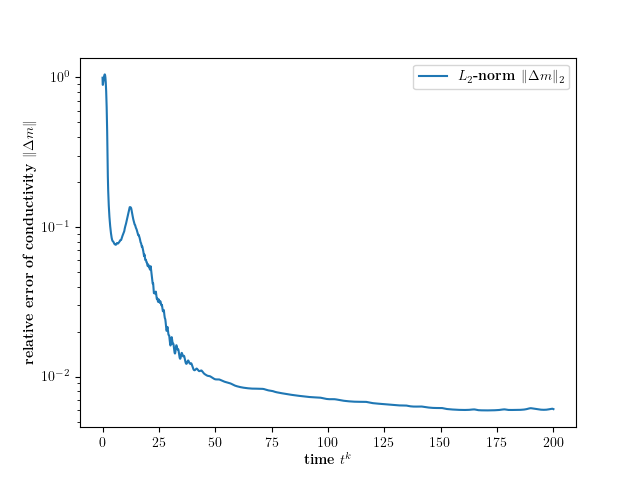}
     \caption{Evolution of $\|\Delta  m\|_2$ over time}
         \label{fig:mini_t3s}
     \end{subfigure}
        \caption[Time evolution of $\cE_h^k$ and $\|\Delta m_h^k$]{Time evolution of the discrete energy $\cE_h^k$ and the relative error of the conductivity vector $\|\Delta m_h^k\|_2$ for $c=200, D=0.0025$, $\gamma=0.75$, $m_{ini}^{const}$ and $S_{int}$ (internal sink). The energy decreases monotonically and stabilises, indicating that the system is close to equilibrium, while the relative change in $m$ is decreasing accordingly.}
        \label{fig:mini_const_energy}
\end{figure}

The corresponding evolution of the discretised energy $\mathcal{E}_h^k$ as well as the evolution of the relative error of the conductivity vector $m$, i.e., $\|\Delta m_h^k\|_2$, over time are shown in Figure \ref{fig:mini_const_energy}. We observe a monotonic decrease in energy and its eventual stabilisation, accompanied by a decrease in the relative change of $m$. Although the breaking condition for the conductivity $m$ is not yet satisfied at $t_{end}=200$, the results suggest that the system is already very close to equilibrium.

\begin{figure}[ht!]
     \centering
     \noindent
\begin{minipage}[t]{0.9\textwidth}
\vspace{0pt}
    \begin{subfigure}[t]{0.32\textwidth}
         \centering
  \includegraphics[width=0.9\linewidth,
                         trim=10mm 5mm 40mm 15mm, clip]{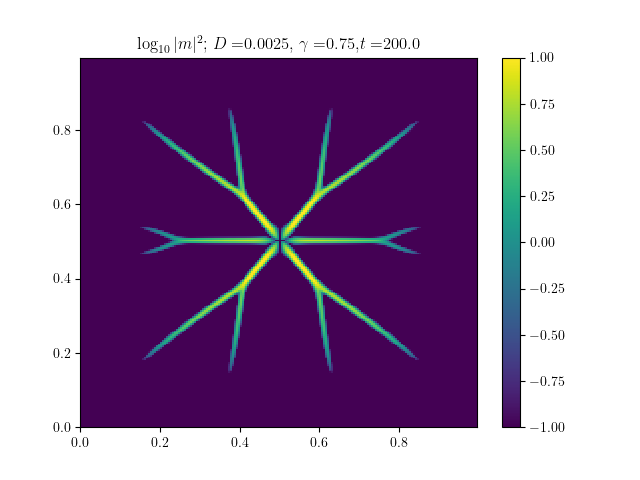}
         \caption{$\log_{10} |m(t,x,y)|$ with $m_{ini}^A$ as initial condition \\}
         \label{fig:mnorm_end_mA}
     \end{subfigure}
     \hfill
     \begin{subfigure}[t]{0.32\textwidth}
         
         \centering
         \includegraphics[width=0.9\linewidth,
                         trim=10mm 5mm 40mm 15mm, clip]      {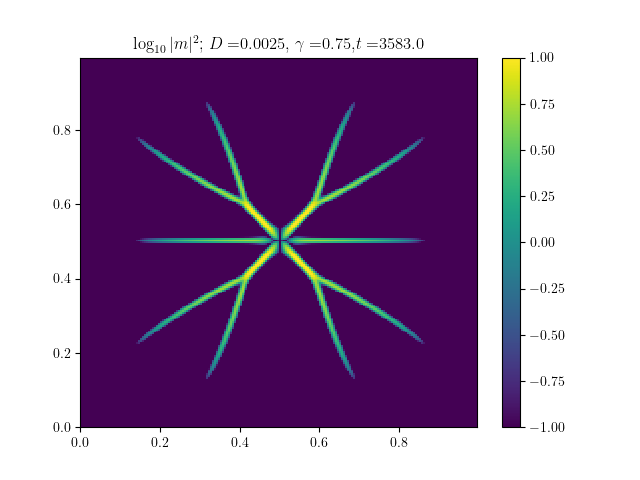}   
     \caption{$\log_{10} |m(t,x,y)|$ with $m_{ini}^B$ as initial condition\\ }
         \label{fig:mnorm_end_mB}
     \end{subfigure}
     \hfill
     \begin{subfigure}[t]{0.32\textwidth}
         \centering
         \includegraphics[width=0.9\linewidth,
                         trim=10mm 5mm 40mm 15mm, clip]      {plots/internal_sink/constant_mini/c200_D00025_gamma075/log_normm_heat_k_40000.png}   
     \caption{$\log_{10} |m(t,x,y)|$ with $m_{ini}^{const}$ as initial condition}
         \label{fig:mnorm_end_mconst}
     \end{subfigure}
     \end{minipage}
     \begin{minipage}[t]{0.05\textwidth} 
     \vspace{-0.2cm}
         \includegraphics[height=3.8cm]{plots/internal_sink/constant_mini/c200_D00025_gamma075/skala.png}
         \text{ } \\
         \text{ }\\
 \end{minipage}
     \noindent
     \begin{minipage}{0.9\textwidth}
         \vspace{0pt}
   
        \begin{subfigure}[t]{0.32\textwidth}
         \centering
         \includegraphics[width=0.9\linewidth,
                         trim=10mm 5mm 40mm 15mm, clip]
         {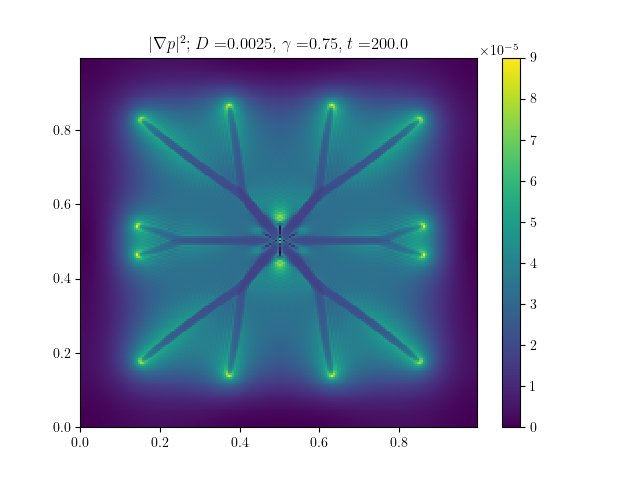}
         \caption{$|\nabla p(t,x,y)|$ with $m_{ini}^A$ as initial condition\\}
         \label{fig:pgrad_norm_end_mA}
     \end{subfigure}
       \hfill
        \begin{subfigure}[t]{0.32\textwidth}
         \centering
         \includegraphics[width=0.9\linewidth,
                         trim=10mm 5mm 40mm 15mm, clip]
         {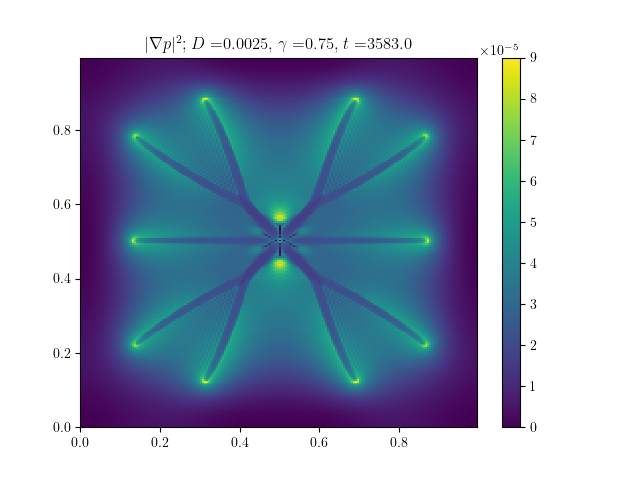}
         \caption{$|\nabla p(t,x,y)|$ with $m_{ini}^B$ as initial condition \\}
         \label{fig:pgrad_norm_end_mB}
     \end{subfigure}
     \hfill
     \begin{subfigure}[t]{0.32\textwidth}
         \centering
         \includegraphics[width=0.9\linewidth,
                         trim=10mm 5mm 40mm 15mm, clip]      {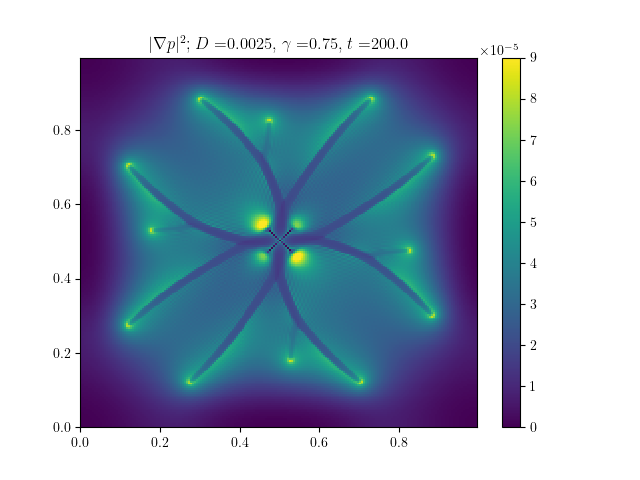}   
     \caption{$|\nabla p(t,x,y)|$ with $m_{ini}^{const}$ as initial condition}
         \label{fig:pgrad_const_end}
     \end{subfigure} 
         \end{minipage}
         \begin{minipage}[t]{0.05\textwidth} 
         \vspace{-3cm}
         \includegraphics[height=3.8cm,
                         trim=125mm 10mm 20mm 5mm, clip]      {plots/internal_sink/constant_mini/c200_D00025_gamma075/gradp_heat_k_40000.png}
 \end{minipage}
       
        \caption[Comparison of simulation results for $m$ and $\nabla p$ for different initial conditions]{Comparison of numerical results for the conductivity $m(t, x, y)$ and the pressure gradient $|\nabla p(t,x,y)|^2$ close to steady state ($t_{end}=200$) for different initial conditions $m_{ini}$. The modelling parameters chosen are $c=200, \gamma= 0.75, r =0.1, D=0.0025$, $\Delta t^k = 0.005$, $S_{int}$ (internal sink). First column: results for $m^A_{ini}$. Middle column: results for $m_{ini}^{B}$. Last column: results for $m_{ini}^{const}$. We conclude that different initial conditions lead to different steady states for $\gamma=0.75$.}
        \label{fig:compare_inicondition}
\end{figure}

Figure \ref{fig:compare_inicondition}
shows a comparison of the final states of the network at $t_{end}=200$ for different initial conditions, i.e., for $m_{ini}^A$ as stated in \eqref{eq:m_ini_normal} (left column of the figure), $m_{ini}^B$ as defined in \eqref{eq:m_ini_symm} (middle column of the figure) and $m_{ini}^{const}$ given in \eqref{eq:m_ini_const} (right column of the figure). It not only depicts the norm of the conductivity vector $|m|$ but also the norm of the pressure gradient $| \nabla p|$ at the final state. Moreover, the different values of the relative error of the conductivity $m$ for the breaking condition are given by $\|\Delta m^{const}\|_2=6.06\cdot 10^{-3}  $ for the constant initial value of $m$ and by  $\|\Delta m^A\|_2=2.67\cdot 10^{-3}$ and $\|\Delta m^B\|_2=3.54\cdot 10^{-3}$ at $t_{end}=200$. Since the energy changes only marginally in all three cases at $t_{end}$, we assume to be close to the steady state in all scenarios. However, the final states for $m$ and $p$ behave very differently depending on the initial conditions. We observe that for the constant initial condition $m_{ini}^{const}$, the branches are curvier, located at different points in space and less clear compared to the results for $m_{ini}^A$ and $m_{ini}^B$. The pressure distribution also differs across all three cases, though several commonalities can be observed: the highest pressure values are consistently located at the tips of the branches and around the source, while the branches themselves exhibit very low pressure. We conclude that different initial conditions might lead to different steady states, which is in coherence with numerical results by \cite{Astuto2022comparison} and the fact that the energy is non-convex.

\begin{figure}[ht!]
     \centering
     \begin{minipage}[t]{0.9\textwidth}
    \begin{subfigure}[b]{0.48\textwidth}
         \centering
         \includegraphics[width=0.9\linewidth,
                         trim=10mm 5mm 40mm 15mm, clip]
         {plots/internal_sink/normal_mini/c200_D00025_gamma075/log_normm_heat_k_40000.png}
         \caption{$\log_{10} |m|$ at $t_{end}=200$ for $c=200$,\linebreak $D=0.0025$}
         \label{fig:mconst_end}
     \end{subfigure}
     \hfill
     \begin{subfigure}[b]{0.48\textwidth}
         \centering
         \includegraphics[width=0.9\linewidth,
                         trim=10mm 5mm 40mm 15mm, clip]      {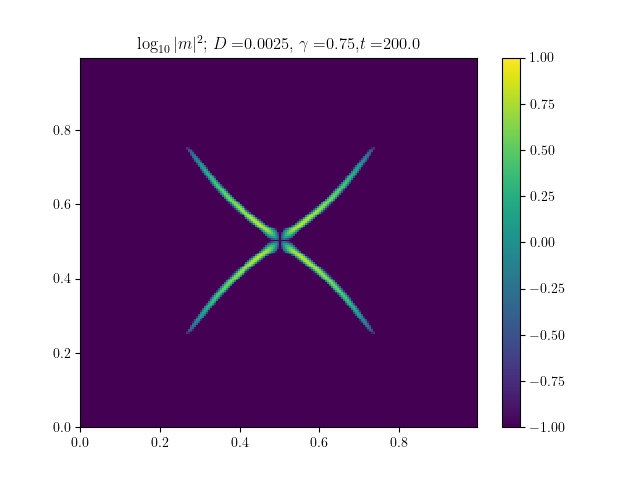}   
     \caption{$\log_{10} |m|$ at $t_{end}=200$ for $c=100$, \linebreak$D=0.0025$}
         \label{fig:mnorm_end}
     \end{subfigure}
     \hfill
        \begin{subfigure}[b]{0.48\textwidth}
         \centering
         \includegraphics[width=0.9\linewidth,
                         trim=10mm 5mm 40mm 15mm, clip]
         {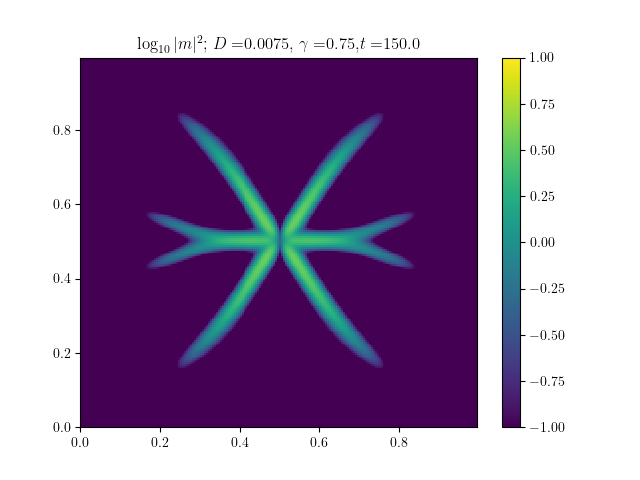}
         \caption{$\log_{10} |m|$ at $t_{end}=150$ for $c=200$, \linebreak $D=0.0075$}
         \label{fig:pgrad_norm_end}
     \end{subfigure}
     \hfill
     \begin{subfigure}[b]{0.48\textwidth}
         \centering
         \includegraphics[width=0.9\linewidth,
                         trim=10mm 5mm 40mm 15mm, clip]      {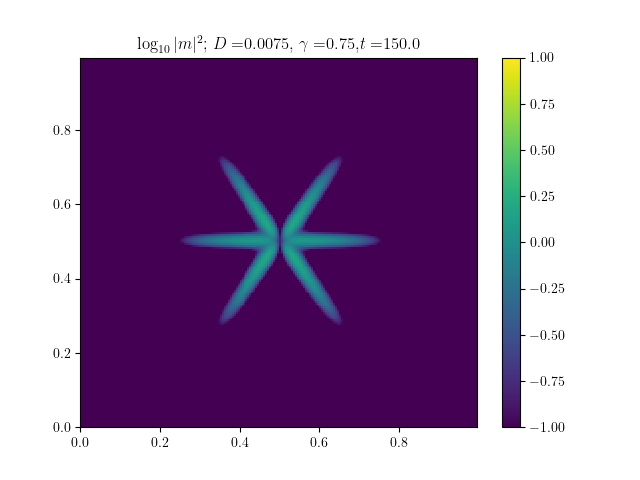}   
     \caption{$\log_{10}|m|$ at $t_{end}=150$ for $c=100$, \linebreak $D=0.0075$}
         \label{fig:lognormm_mA_end}
     \end{subfigure} 
      \end{minipage}
      \begin{minipage}{0.05\textwidth}   
       \begin{subfigure}{0.9\textwidth}
         \centering
         \includegraphics[height=8.8cm]
         {plots/internal_sink/constant_mini/c200_D00025_gamma075/skala.png}
     \end{subfigure}
 \end{minipage}
        \caption[Comparison of simulation results of $m$ for varying activation parameter $c$]{Comparison of numerical results for the conductivity $m$ close to steady state for different values of the activation parameter $c$. The parameters chosen are $\gamma=0.75, r=0.1, \Delta t^k=0.005$, $m_{ini}^A$ and the internal sink $S_{int}$. Left column: $c=200$. Right column: $c=100$. Top row: $D=0.0025$. Bottom row: $D=0.0075$. Larger activation parameters result in thicker, more branched structures that extend further into the domain.}
        \label{fig:compare_cvalues}
\end{figure}

Next, we investigate the influence of the activation parameter $c$ on the network. Therefore, we consider the values $c=200$ and $c=100$ for two different setups. In both scenarios, we use $\gamma=0.75$ as the value for the exponent of the metabolic cost term. On the other hand, we consider $D=0.0025$ and $D=0.0075$ as the diffusion constant. The numerical results for the logarithm of the Euclidean norm of the conductivity $m(t,x,y)$ are depicted in Figure \ref{fig:compare_cvalues}. The right column of the figure shows the results for $c=100$, and the left column shows $ c=200$. The first row corresponds to $D=0.0025$, and the second row corresponds to $D=0.0075$. We observe that for $c=200$ the network forms thicker and more pronounced branches, which also reach further into the domain than for $c=100$. On the other hand, smaller activation ($c=100$) leads to sparser networks, and in the presented cases, the networks even collapse to a simple cross (for $D=0.0025$) or star-shaped configuration (for $D=0.0075$). Moreover, we conclude that for smaller diffusion parameters, as depicted in the top row for $D=0.0025$, the networks are finer and more detailed with thin filaments extending outward, while larger diffusion ($D=0.0075$, bottom row) leads to smoother and more compact networks with fewer but thicker branches. Hence, the combined effects indicate that activation has to dominate the diffusion part for the emergence of complex networks, while increased diffusion or reduced activation suppresses branching and the system tends to converge to simpler steady states.  \\

\begin{figure}[ht!]
     \centering
\hspace{-1.5cm}
     \begin{minipage}{0.9\textwidth}  
     \hfill
    \begin{subfigure}[b]{0.18\textwidth}
         \centering
         \includegraphics[width=0.9\linewidth,
                         trim=10mm 0mm 40mm 15mm, clip]     {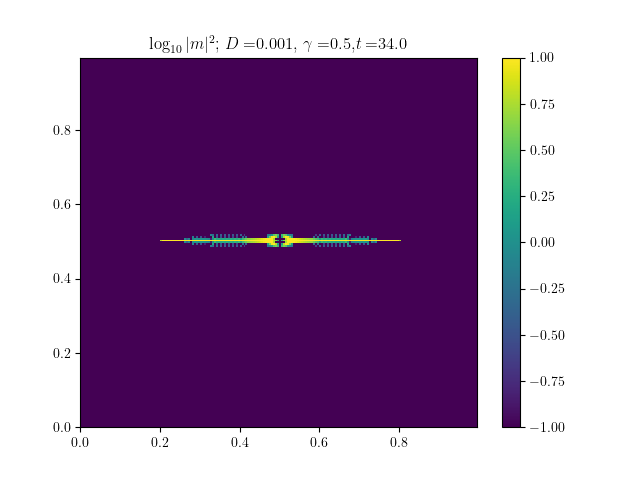}
     \end{subfigure}
     \hfill
     \begin{subfigure}[b]{0.18\textwidth}
         \centering
         \includegraphics[width=0.9\linewidth,
                         trim=10mm 0mm 40mm 15mm, clip]         {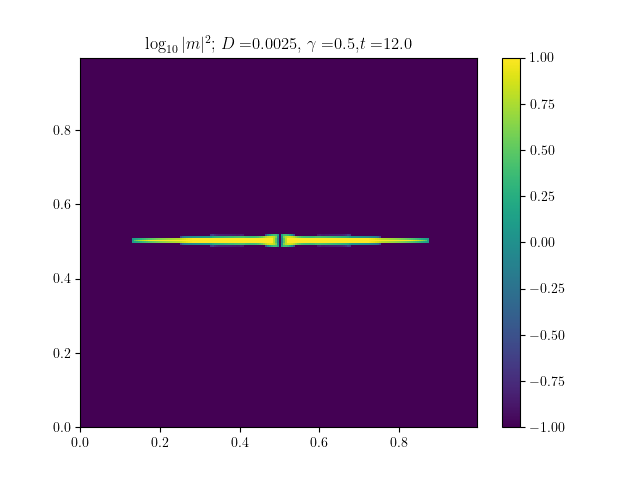}
     \end{subfigure}
      \hfill
     \begin{subfigure}[b]{0.18\textwidth}
         \centering
         \includegraphics[width=0.9\linewidth,
                         trim=10mm 0mm 40mm 15mm, clip]        {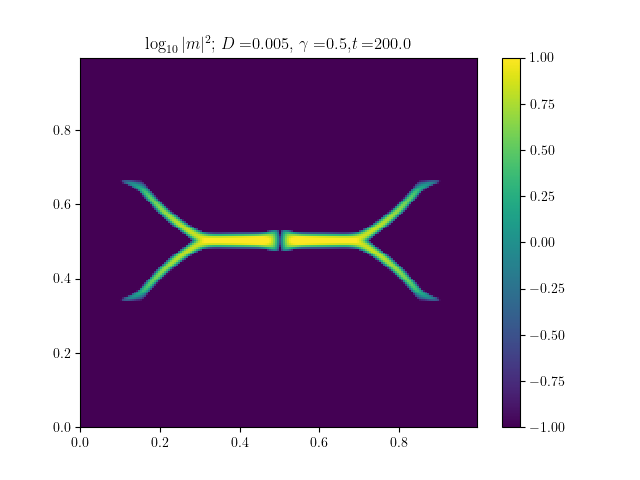}
     \end{subfigure}
      \hfill
     \begin{subfigure}[b]{0.18\textwidth}
         \centering
         \includegraphics[width=0.9\linewidth,
                         trim=10mm 0mm 40mm 15mm, clip]         {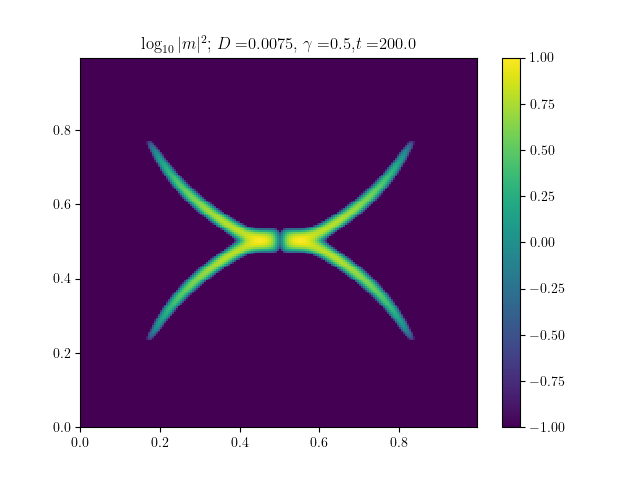}
     \end{subfigure}
     \hfill
     \begin{subfigure}[b]{0.18\textwidth}
         \centering
         \includegraphics[width=0.9\linewidth,
                         trim=10mm 0mm 40mm 15mm, clip]
         {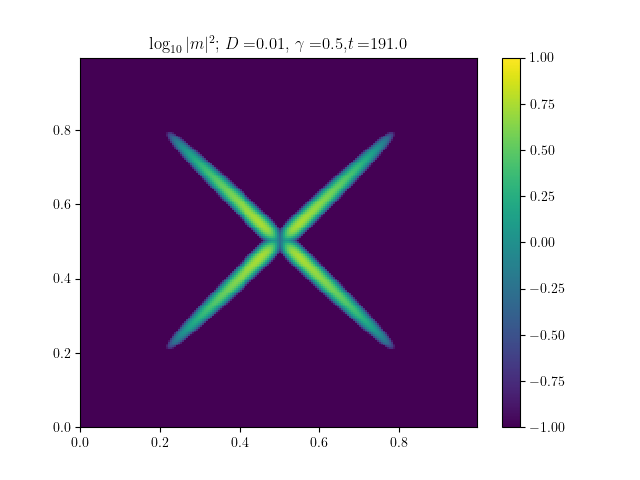}
     \end{subfigure}

     \hfill
        \begin{subfigure}[b]{0.18\textwidth}
         \centering
         \includegraphics[width=0.9\linewidth,
                         trim=10mm 0mm 40mm 15mm, clip]        {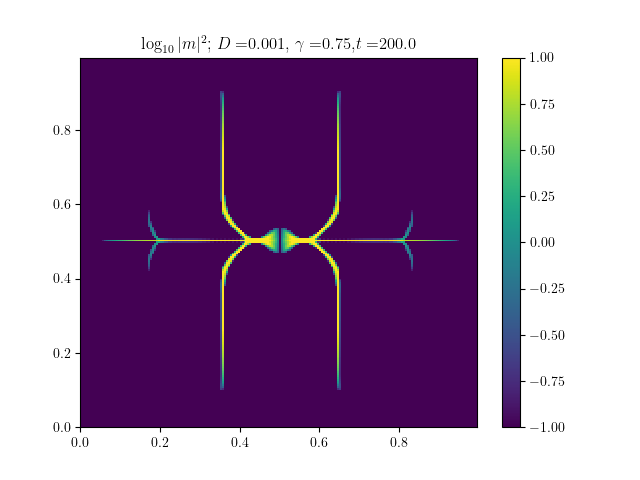}
     \end{subfigure}
     \hfill
     \begin{subfigure}[b]{0.18\textwidth}
         \centering
         \includegraphics[width=0.9\linewidth,
                         trim=10mm 0mm 40mm 15mm, clip]         {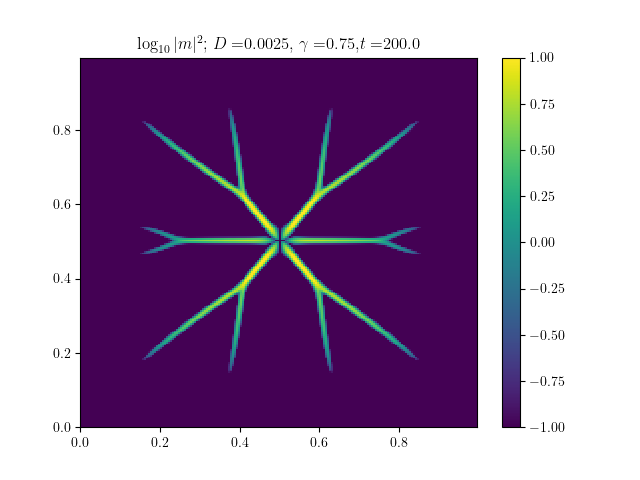}
     \end{subfigure}
      \hfill
     \begin{subfigure}[b]{0.18\textwidth}
         \centering
         \includegraphics[width=0.9\linewidth,
                         trim=10mm 0mm 40mm 15mm, clip]        {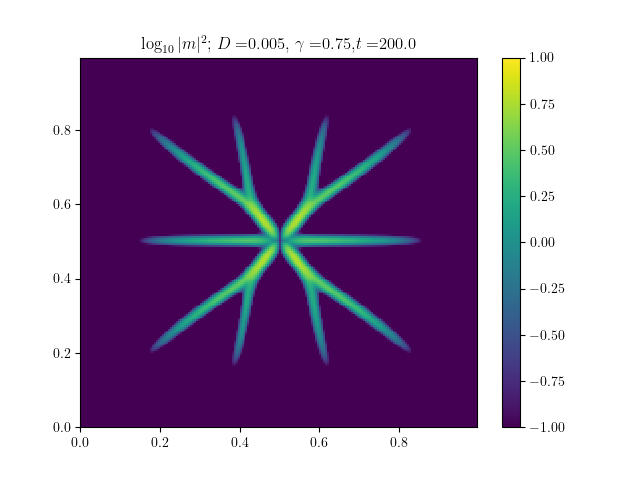}
     \end{subfigure}
      \hfill
     \begin{subfigure}[b]{0.18\textwidth}
         \centering
         \includegraphics[width=0.9\linewidth,
                         trim=10mm 0mm 40mm 15mm, clip]         {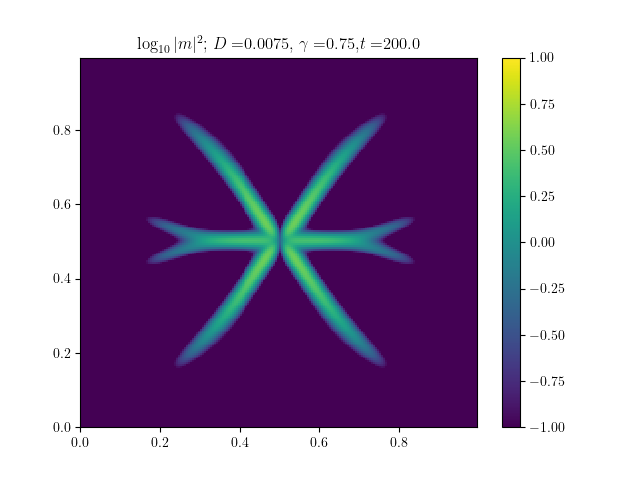}
     \end{subfigure}
     \hfill
     \begin{subfigure}[b]{0.18\textwidth}
         \centering
         \includegraphics[width=0.9\linewidth,
                         trim=10mm 0mm 40mm 15mm, clip]
         {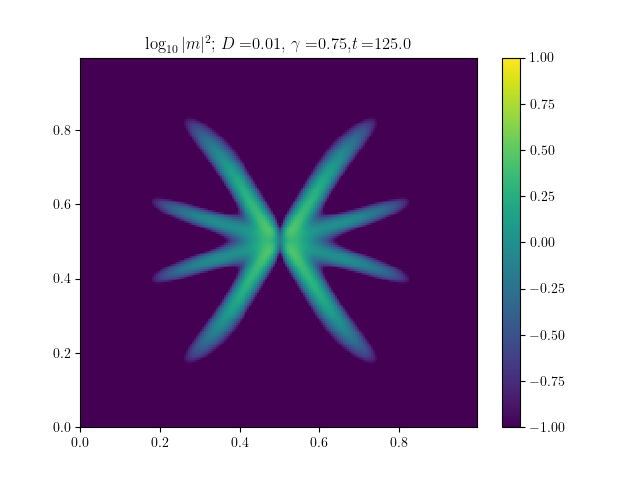}
     \end{subfigure}

      \hfill
        \begin{subfigure}[b]{0.18\textwidth}
         \centering
         \includegraphics[width=0.9\linewidth,
                         trim=10mm 0mm 40mm 15mm, clip]        {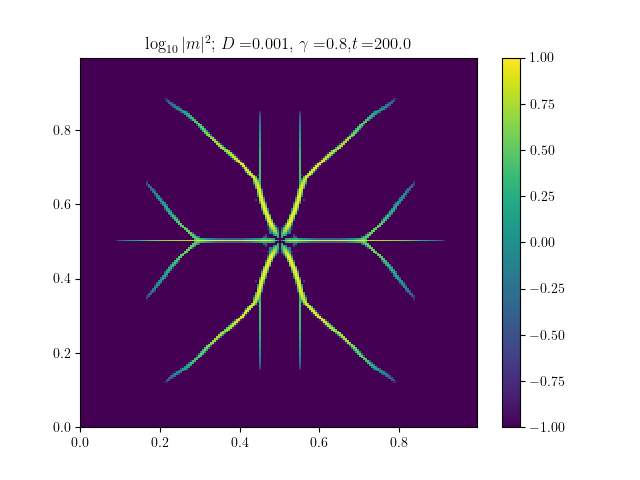}
     \end{subfigure}
     \hfill
     \begin{subfigure}[b]{0.18\textwidth}
         \centering
         \includegraphics[width=0.9\linewidth,
                         trim=10mm 0mm 40mm 15mm, clip]         {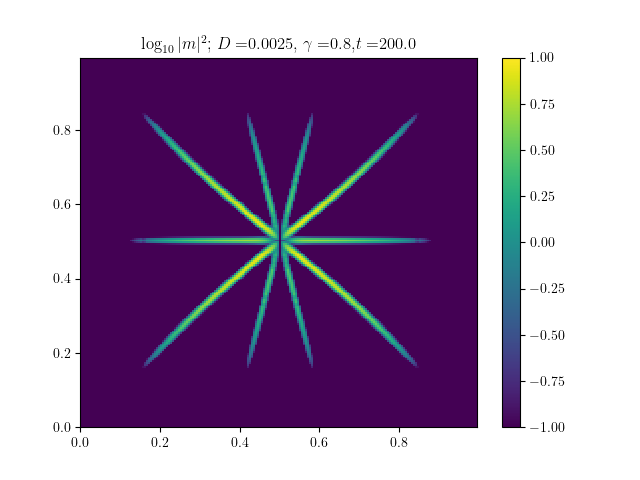}
     \end{subfigure}
      \hfill
     \begin{subfigure}[b]{0.18\textwidth}
         \centering
         \includegraphics[width=0.9\linewidth,
                         trim=10mm 0mm 40mm 15mm, clip]       {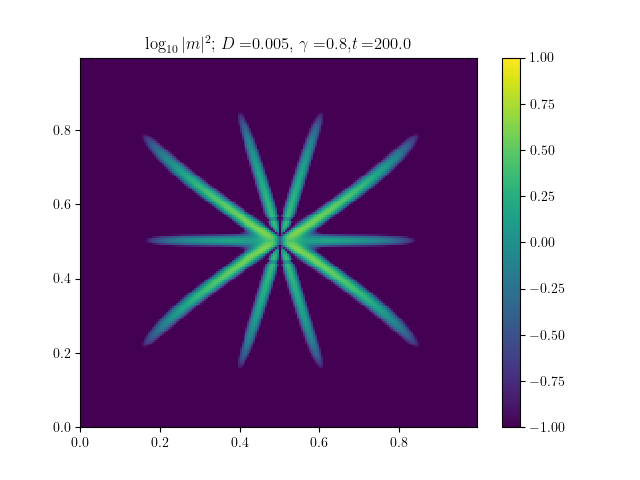}
     \end{subfigure}
      \hfill
     \begin{subfigure}[b]{0.18\textwidth}
         \centering
         \includegraphics[width=0.9\linewidth,
                         trim=10mm 0mm 40mm 15mm, clip]        {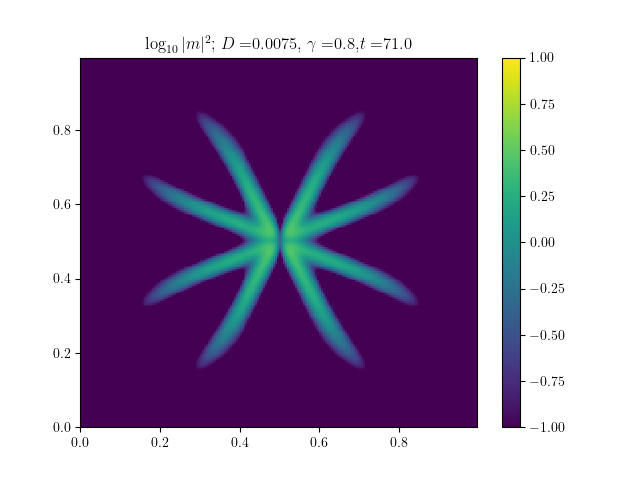}
     \end{subfigure}
     \hfill
     \begin{subfigure}[b]{0.18\textwidth}
         \centering
         \includegraphics[width=0.9\linewidth,
                         trim=10mm 0mm 40mm 15mm, clip]
         {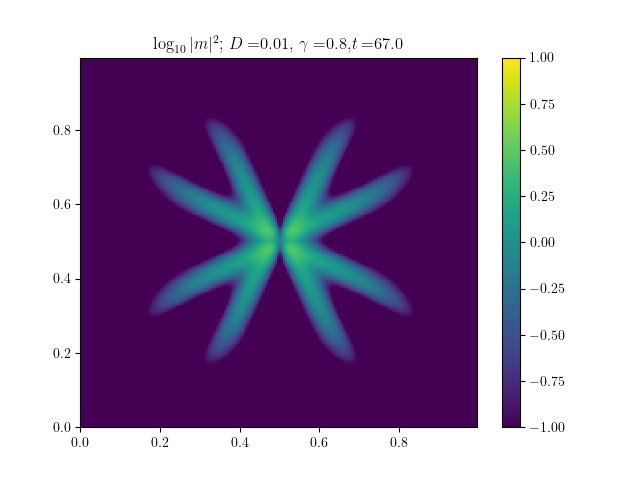}
     \end{subfigure}

         \hfill
        \begin{subfigure}[b]{0.18\textwidth}
         \centering
         \includegraphics[width=0.9\linewidth,
                         trim=10mm 0mm 40mm 15mm, clip]         {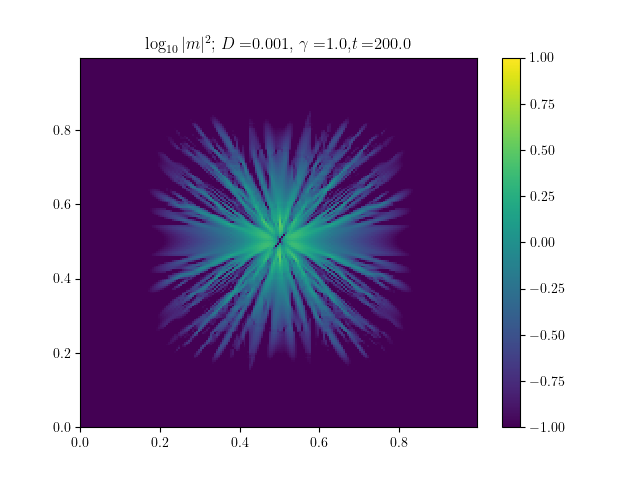}
     \end{subfigure}
     \hfill
     \begin{subfigure}[b]{0.18\textwidth}
         \centering
         \includegraphics[width=0.9\linewidth,
                         trim=10mm 0mm 40mm 15mm, clip]       {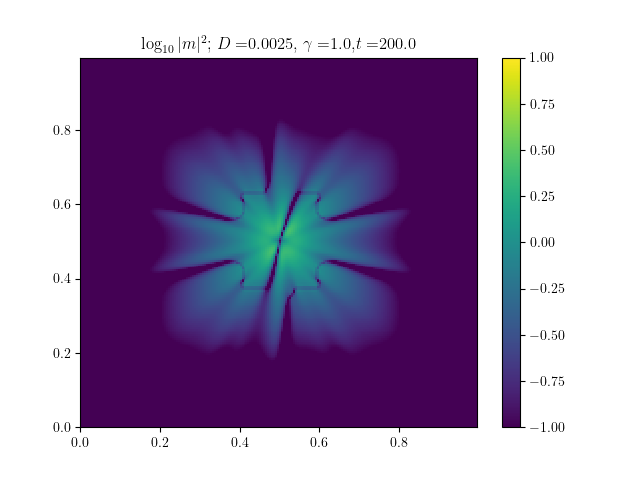}
     \end{subfigure}
      \hfill
     \begin{subfigure}[b]{0.18\textwidth}
         \centering
         \includegraphics[width=0.9\linewidth,
                         trim=10mm 0mm 40mm 15mm, clip]        {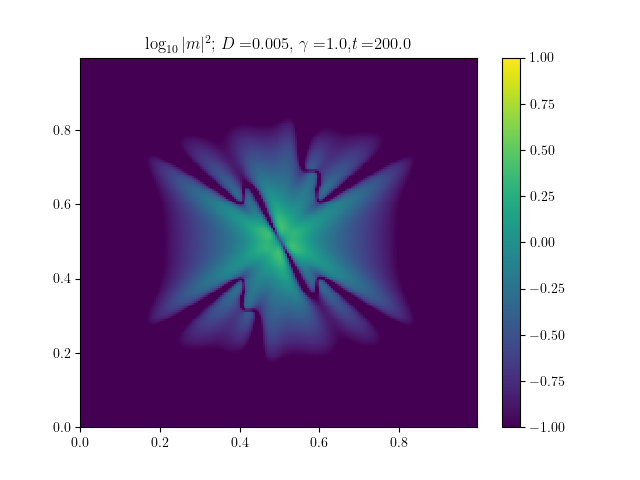}
     \end{subfigure}
      \hfill
     \begin{subfigure}[b]{0.18\textwidth}
         \centering
         \includegraphics[width=0.9\linewidth,
                         trim=10mm 0mm 40mm 15mm, clip]        {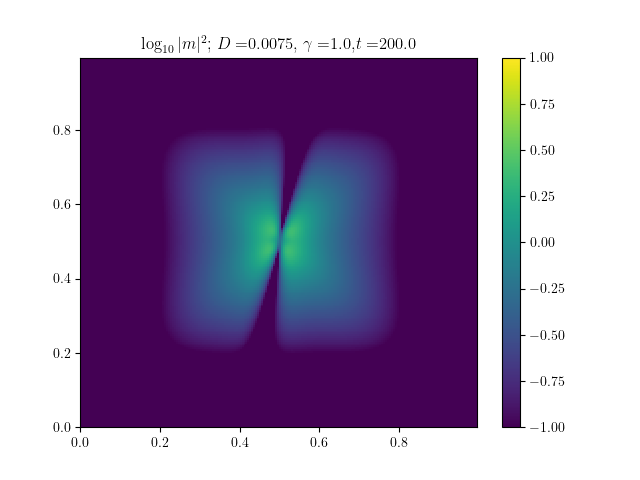}
     \end{subfigure}
     \hfill
     \begin{subfigure}[b]{0.18\textwidth}
         \centering
         \includegraphics[width=0.9\linewidth,
                         trim=10mm 0mm 40mm 15mm, clip]
         {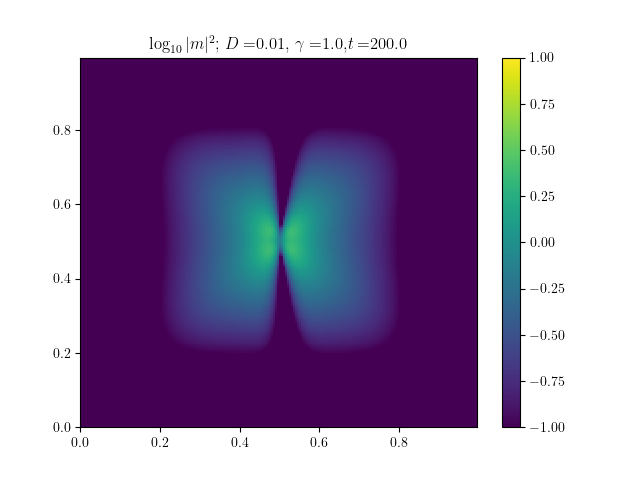}
     \end{subfigure}
     \end{minipage}
\begin{minipage}{0.05\textwidth}
         \begin{subfigure}[t]{1.5\textwidth}
         \centering
         \vspace{-0.2cm}
         \includegraphics[height=9.5cm]
         {plots/internal_sink/constant_mini/c200_D00025_gamma075/skala.png}
     \end{subfigure}
     \end{minipage}     
        \caption[Parameter study for $\gamma$ and $D$]{Value $\log_{10} |m|$ for $c=200$, $m^A_{ini}$, and $S_{int}$ (internal sink). From top to row $\gamma$ varies in $\{0.5,0.75,0.8,1\}$ and from left to right columns $D$ varies in $\{0.001, 0.0025,0.005,0.0075,0.01\}$. The network shown is the one at time $t_{end}=200$ or when the algorithm broke due to an increase in energy or reaching the steady state.}
        \label{fig:compare_D_gamma}
\end{figure}

In Figure \ref{fig:compare_D_gamma} we study the influence of the diffusion parameter $D$ and the metabolic exponent $\gamma$ on the evolution of the network. We consider $c = 200$ as the activation constant. The initial condition for the conductivity $m(t,x,y)$ is given by $m_{ini}^A$ as stated in \eqref{eq:m_ini_normal}. The numerical results for $\log_{10}|m|$ of the conductivity $m$ are depicted in that figure, where $\gamma$ varies in $\{0.5,0.75,0.8,1\}$ (top to bottom row of the figure) and $D\in \{0.001, 0.0025,0.005,0.0075,0.01\}$ (left to right column of the figure). We observe that
for small values of diffusion, finer and highly branched networks evolve. On the other hand, increasing $D$ leads to thicker and smoother structures, suppressing the fine details of the network. Hence, smaller diffusion constants form rich branching patterns, while larger diffusion enforces global uniformity. Moreover, we conclude from the figure that the exponent of the metabolic cost term $\gamma$ controls the density of the network, that is, small values of $\gamma$ give rise to sparse and tree-like patterns since the system minimises total metabolic cost, while increasing the value of $\gamma$ leads to denser and more space-filling networks, and thus to broader regions of active flow. Overall, we deduce from these simulation results that complex networks emerge when $ \gamma$ has a moderate value (i.e., $\gamma=0.75$) and small values for the diffusion coefficient $D$. Large diffusion (i.e., $D=0.01$) and large $\gamma$ lead to nearly homogeneous steady states.   

\section{Grid Convergence}
\label{sec:gridconv}

To investigate whether the simulation outcomes depend on the chosen grid, as observed in \cite{Haskovec2016}, we perform a grid convergence test. For this, we will take values of $c$ corresponding to $c=100$ and $c=50$ to reduce the computational time. We compute the order of convergence for the pressure $p$ via the following formula
\begin{equation}
	q_p\approx \frac{\ln \left( \frac{\|p_{\Delta x}-p_{\frac{\Delta x}{2}}\|_2 }{\|p_{\frac{\Delta x}{2}}-p_{\frac{\Delta x}{4}}\|_2}\right)}{\ln 2},
\end{equation}
where $p_{\Delta x}$ denotes the solution for $p$ on a regular grid with grid size $\Delta x$.
We define analogously $q_{\|m\|}$ and $q_{\|\nabla p\|}$.

We consider $\Omega=[0,1]^2$, $N_x=N_y \in\{50,100,200,400\}$, $r=0.1$, and $S=S_{\mathrm{int}}$. Table \ref{tab:q_gamma075} corresponds to the setup with  $c=100$, $\gamma=0.75$, $D=0.01$, $\Delta t=0.05$ and  $m^A_{\mathrm{ini}}$ and reports the rate $q$ at increasing times  (the coarsest run terminates at $t\approx 91$ due to reaching steady state). Table \ref{tab:q_gamma2} considers the case $c=50$, $\gamma=1$, $D=0.0025$ and $m^{const}_{ini}$. Table \ref{tab:q_gamma3} considers the case  with $m_{ini}^{A}$ with factor $1000$ instead of $10^4$ in the exponential, $c=50, \gamma=1, D=0.0025$. Missing numbers are either due to the simulations still running (i.e., Table \ref{tab:q_gamma075}) 
or that they terminated at an earlier time (i.e., Table \ref{tab:q_gamma2} for $N=50,100,200$ and $t=35$ and $t=40$ as well as $t=40$ for $N=100,200,400$).

 \begin{table}[ht!]
 	\centering
 	\begin{tabular}{|c|c|c|c||c|c|c||c|c|c|}
 		\hline
 		& \multicolumn{3}{|c||}{$N=50,100,200$} &\multicolumn{3}{|c|}{$N=100,200,400$} &\multicolumn{3}{|c|}{$N=200,400,800$} \\
 		\hline
 		time $t$  & $q_{\| m\|}$ & $q_{\| \nabla p \|}$ & $q_p$ & $q_{\| m\|}$ & $q_{\| \nabla p \|}$ & $q_p$ & $q_{\| m\|}$ & $q_{\| \nabla p \|}$ & $q_p$\\
 		\hline
 		\hline
        5 & 6.6747 & 5.0983 & 5.6516 & 3.9922&  3.8377 & 3.2919 & 3.4010 & 4.5286 & 5.6811\\
 		10 & 6.9849 & 5.0677& 5.6942 & 3.7730  & 3.9284 &  3.4008 &--- &--- &--- \\
 		50 & 7.9296 & 6.3999 & 7.8724 & 3.6918 & 3.6324  & 2.8776 & ---& --- & --- \\
        	91 & 7.8778& 6.3763 & 7.6062& 4.2736 & 4.1772  & 4.0374 & --- & --- & ---\\		
 		\hline
 		\hline 
 	\end{tabular} 			
 	\caption{Observed orders $q$ for $c=100$, $\gamma=0.75$, $D=0.01$, $\Delta t=0.05$, $m^A_{\mathrm{ini}}$,  $S_{\mathrm{int}}$, $N=N_x=N_y.$}
 		\label{tab:q_gamma075}
 	\end{table}

 \begin{table}[ht!]
 	\centering
 	\begin{tabular}{|c|c|c|c||c|c|c||c|c|c|}
 		\hline
    & \multicolumn{3}{|c||}{$N=50,100,200$} &\multicolumn{3}{|c||}{$N=100,200,400$} &\multicolumn{3}{|c|}{$N=200,400,800$}\\
 		\hline
 		time $t[s]$  & $q_{\| m\|}$ & $q_{\| \nabla p \|}$ & $q_p$ & $q_{\| m\|}$ & $q_{\| \nabla p \|}$ & $q_p$ & $q_{\| m\|}$ & $q_{\| \nabla p \|}$ & $q_p$\\
 		\hline 
 		\hline
 		15 & 2.5533 & 2.209 & 2.7430 & 4.6220 & 3.7932  & 4.4483 & 8.2958 &6.2632 & 6.8885  \\
        	35 &---& --- & --- & 4.7246 & 3.7308 & 4.5646& 8.3799 & 6.4396 & 6.9787 \\
            40 & --- &---&---& ---&---&---&8.3838 & 6.4438 & 6.9735\\
 		\hline
 		\hline 
 	\end{tabular} 			
 	\caption{Observed orders $q$ for $c=100,\ \gamma=1,\ D=0.0025,\ \Delta t=0.5$, $m_{ini}^{const}$, $S_{\mathrm{int}}$, $N=N_x=N_y.$}
 		\label{tab:q_gamma2}
 	\end{table}

 \begin{table}[ht!]
 	\centering
 	\begin{tabular}{|c|c|c|c||c|c|c||c|c|c|}
 		\hline
    & \multicolumn{3}{|c||}{$N=50,100,200$} &\multicolumn{3}{|c||}{$N=100,200,400$} &\multicolumn{3}{|c|}{$N=200,400,800$}\\
 		\hline
 		time $t[s]$  & $q_{\| m\|}$ & $q_{\| \nabla p \|}$ & $q_p$ & $q_{\| m\|}$ & $q_{\| \nabla p \|}$ & $q_p$ & $q_{\| m\|}$ & $q_{\| \nabla p \|}$ & $q_p$\\
 		\hline 
 		\hline
        45 &4.8019& 2.5699 &3.773& 2.5099& 2.7235 & 2.3999& 6.6369 & 3.9876 & 5.5382 \\
 		55 & 4.8009 & 2.5625 & 3.7567 & 2.5115 & 2.7251  & 2.3969 & 6.639 & 3.987 & 5.5368  \\
 		\hline
 		\hline 
 	\end{tabular} 			
 	\caption{Observed orders $q$ for $c=100,\ \gamma=1,\ D=0.0025,\ \Delta t=0.05$ with $m_{ini}^A$ with factor $1000$ instead of $10^4$, $S_{\mathrm{int}}$, $N=N_x=N_y.$}
 		\label{tab:q_gamma3}
 	\end{table}

The numerical method proposed utilises the FFT, the CG method, and the splitting method. The splitting method has order 1, whereas in the CG method, we can choose the desired accuracy. Spectral methods have an exponential order of convergence. Due to this mixture of accuracies, it is not obvious to predict which should be the expected order of convergence of the method. However, so far, the numerical evidence shows that the method has a high order of convergence. It is in particular noticeable that no observable aliasing effects appear in the simulations.

\section{Conclusions and Perspectives}\label{sec:conclusion_network_formation}
In this paper, we considered the elliptic-parabolic model \eqref{eq:networkSystem} to study the formation of biological networks. The model was derived as a continuum limit of the discrete Hu-Cai model and formulated as an $L^2$-gradient flow of a non-convex energy functional.

We introduced a Fourier-based spectral splitting method for the network formation PDE \eqref{eq:networkSystem} under periodic boundary conditions. The algorithm is built from (1) a splitting method for the conductivity $m$ with
(i) exact-in-time substeps for activation and relaxation,
(ii) an exact diffusion step implemented with FFT, and (2) a CG method to solve the equation for the pressure $p$ based on using FFT. 

This combination yields a simple-to-implement solver with high spatial accuracy and scalability to 3D uniform grids. Our parameter studies recover qualitative features reported in the literature: diffusion $D$ smooths and thickens networks, activation $c$ enhances branching, and the metabolic exponent $\gamma$ tunes the sparsity/density of steady patterns. Grid convergence results support the accuracy of our approach. However, the method can be improved in different directions, in particular, it is not unconditionally energy-stable and it is currently limited by 
a time-step restriction due to the activation update (Remark \ref{rem:Delta_t}), for which adaptive time-stepping strategies may be used in practice. Also, to reduce the computational time, preconditioners beyond the constant-coefficient could further reduce CG iterations. 

Despite its limitations, the proposed method fills a practical niche among existing solvers: it is concise, fast, and robust enough for both qualitative and quantitative studies, and is easily extensible to 3D problems where network formation is particularly relevant for applications such as vascularisation.

\section{Acknowledgments}

The authors would like to thank Steffen Plunder for wonderful and fruitful discussions.
\medskip

The work of SMA and CW was funded in part by the Austrian Science Fund (FWF) project \href{https://www.fwf.ac.at/en/research-radar/10.55776/F65}{10.55776/F65}
and in part by the Vienna Science and Technology Fund (WWTF) \href{https://www.wwtf.at/funding/programmes/vrg/VRG17-014/index.php?lang=EN}{[10.47379/VRG17014]}. CW was also funded partly by the Austrian Science Fund (FWF) project \href{https://www.fwf.ac.at/forschungsradar/10.55776/W1261}{10.55776/W1261}. PAS acknowledges support from the Simons Foundation through the Travel Support for Mathematicians program.

\appendix

\section{Videos of Simulations}

\label{sec:videos}
\textbf{Videos 1-4: Evolution of the conductance vector $m$ and the pressure gradient $p$ for the setup with four sinks and one source} 

\textit{Description:} These videos show the evolution of the network via the conductance $m$ and the pressure gradient $p$ for the parameters $c=200, D=0.0025$ and $\gamma=0.7$ as well as $\gamma=0.1$ and for the initial setups $m_{ini}^A$ and $S_4$ (one source, four sinks). These simulations highlight the transition of branching as well as observing more symmetric networks for larger $\gamma$ and supplements Figures \ref{fig:4sinks_g07}, \ref{fig:4sinks_g1}, \ref{fig:4sinks_g07_gradp} and \ref{fig:4sinks_g1_gradp}.\\

\noindent\textbf{Video 5: Evolution of the conductance vector $m$ for Figure \ref{fig:mini_time_evolution}}

\textit{Description:} The video shows the network formation over time via the conductance $m$ for the parameters $c=200, D=0.0025$ and $\gamma=0.75$ and the initial setup of the internal sink $S_{int}$ and the initial constant conductivity $m_{ini}^{const}$. This video supplements Figure \ref{fig:mini_time_evolution}.\\

\noindent\textbf{Videos 6-8: Evolution of the conductance vector $m$ three different cases of Figure \ref{fig:compare_D_gamma}}

\textit{Description:} The videos show the different time evolution of the network via the conductance $m$ for the fixed parameter $c=200$ and the varying parameters $(D=0.0025, \gamma=0.75)$, $(D=0.005,\gamma=0.8)$ and $(D=0.0075, \gamma=0.75)$ in the case of the internal sink $S_{int}$ and the initial condition for $m$ given by $m_{ini}^A$. The videos highlight the influence of the parameters $D$ and $\gamma$ on the emergence of the network and support Figure \ref{fig:compare_D_gamma}.\\

For all videos see link: \url{https://doi.org/10.6084/m9.figshare.32433996}

\bibliographystyle{plain}
\bibliography{main}

\end{document}